\documentclass[12pt,english]{amsart}
\usepackage[T1]{fontenc}
\usepackage[latin9]{inputenc}
\usepackage{geometry}
\usepackage{color}
\usepackage{babel}
\usepackage{verbatim}
\usepackage{amsthm}
\usepackage{amstext}
\usepackage{amssymb}
\usepackage{esint}
\usepackage[unicode=true,pdfusetitle,
 bookmarks=true,bookmarksnumbered=false,bookmarksopen=false,
 breaklinks=false,pdfborder={0 0 1},backref=section,colorlinks=false]
 {hyperref}
\usepackage{breakurl}

\makeatletter
\numberwithin{equation}{section}
\theoremstyle{plain}
\newtheorem{thm}{\protect\theoremname}[section]
  \theoremstyle{definition}
  \newtheorem{defn}[thm]{\protect\definitionname}
  \theoremstyle{plain}
  \newtheorem{fact}[thm]{\protect\factname}
  \theoremstyle{remark}
  \newtheorem{rem}[thm]{\protect\remarkname}
  \theoremstyle{remark}
  \newtheorem*{rem*}{\protect\remarkname}
  \theoremstyle{plain}
  \newtheorem{lem}[thm]{\protect\lemmaname}
  \theoremstyle{remark}
  \newtheorem{claim}[thm]{\protect\claimname}
  \theoremstyle{plain}
  \newtheorem{prop}[thm]{\protect\propositionname}

\usepackage{babel}

  \providecommand{\definitionname}{Definition}
  \providecommand{\factname}{Fact}
  \providecommand{\lemmaname}{Lemma}
  \providecommand{\propositionname}{Proposition}
  \providecommand{\remarkname}{Remark}
\providecommand{\theoremname}{Theorem}

\makeatother

  \providecommand{\claimname}{Claim}
  \providecommand{\definitionname}{Definition}
  \providecommand{\factname}{Fact}
  \providecommand{\lemmaname}{Lemma}
  \providecommand{\propositionname}{Proposition}
  \providecommand{\remarkname}{Remark}
\providecommand{\theoremname}{Theorem}

\begin{document}

\title[Interpolation spaces between $\ell^{1}$ and $\ell^{q}$ ]{An alternative characterization of normed interpolation spaces between
$\ell^{1}$ and $\ell^{q}$}

\author{Michael Cwikel}

\address{Cwikel: Department of Mathematics, Technion - Israel Institute of
Technology, Haifa 32000, Israel }

\email{mcwikel@math.technion.ac.il }

\author{Per G. Nilsson}

\address{Nilsson : Wollmar Yxkullsgatan 5A, 118 50 Stockholm, Sweden }

\email{pgn@plntx.com}

\thanks{The first named author's work was supported by the Technion V.P.R.\ Fund
and by the Fund for Promotion of Research at the Technion.}


\begin{abstract}
Given a constant $q\in(1,\infty)$, we study the following property
of a normed sequence space $E$: 

\textit{If $\left\{ x_{n}\right\} _{n\in\mathbb{N}}$ is an element
of $E$ and if $\left\{ y_{n}\right\} _{n\in\mathbb{N}}$ is an element
of $\ell^{q}$ such that
\[
\sum_{n=1}^{\infty}\left|x_{n}\right|^{q}=\sum_{n=1}^{\text{\ensuremath{\infty}}}\left|y_{n}\right|^{q}
\]
and if the nonincreasing rearrangements of these two sequences satisfy
\[
\sum_{n=1}^{N}\left|x_{n}^{*}\right|^{q}\le\sum_{n=1}^{N}\left|y_{n}^{*}\right|^{q}
\]
for all $N\in\mathbb{N}$, then $\left\{ y_{n}\right\} _{n\in\mathbb{N}}\in E$
and $\left\Vert \left\{ y_{n}\right\} _{n\in\mathbb{N}}\right\Vert _{E}\le C\left\Vert \left\{ x_{n}\right\} _{n\in\mathbb{N}}\right\Vert _{E}$
for some constant $C$ which depends only on $E$.}

We show that this property is very close to characterizing the normed
interpolation spaces between $\ell^{1}$ and $\ell^{q}$. More specificially,
we first show that every space which is a normed interpolation space
with respect to the couple $\left(\ell^{p},\ell^{q}\right)$ for some
$p\in[1,q]$ has the above mentioned property. Then we show, conversely,
that if $E$ has the above mentioned property, and also has the Fatou
property, and is contained in $\ell^{q}$, then it is a normed interpolation
space with respect to the couple $\left(\ell^{1},\ell^{q}\right)$.
These results are our response to a conjecture of Galina Levitina,
Fedor Sukochev and Dmitriy Zanin in \cite{LeSuZa2017B}.
\end{abstract}
\maketitle

\section{Acknowledgements}

This research was begun while both the authors were participants in
the conference on Banach Spaces and Operator Theory with Applications,
on the occasion of the 60th birthday of Professor Mieczys{\l}aw Masty{\l}o,
at Adam Mickiewicz University in Pozna\'n, Poland. We are very grateful
to Adam Mickiewicz University and to the organizers of the conference,
in particular Professor Pawe{\l} Mleczko, for their fine hospitality
and creation of an excellent environment for discussions and research. 

We are extremely grateful to Professor Fedor Sukochev who motivated
us with his very interesting lecture at the conference, in which he
discussed topics in his joint work \cite{LeSuZa2017,LeSuZa2017B}
with Galina Levitina and Dmitriy Zanin and, in particular, made us
aware of a remarkable conjecture of theirs. This paper is our response
to that conjecture. We also thank him for some helpful subsequent
discussions.

\section{\label{sec:Intro}Introduction and some auxiliary results}

In most cases where it is has been found to be possible to describe
all the interpolation spaces with respect to a given Banach couple,
the description is formulated via a monotonicity condition for the
$K$-functional for that couple. The terminology \textit{Calder\textit{\'o}n
couple} or \textit{Calder\textit{\'o}n-Mityagin couple} is often
used to describe those couples for which such a description of all
their interpolation spaces is valid.

In this paper we consider a couple which is already known to be a
Calder\'on couple, namely $(\ell^{1},\ell^{q})$ where $q\in(1,\infty)$.
We characterize all the interpolation spaces of this couple, or at
least those which are normed and have the Fatou property, via the
following different property which we refer to as $S_{q}(C)$, and
which is expressed in terms of a functional which is essentially different
from the $K$-functional for $\left(\ell^{1},\ell^{q}\right)$, but,
as we shall see later, does have some kind of connection with that
$K$-functional. 
\begin{defn}
\label{def:SqC} Let $E\ne\left\{ 0\right\} $ be a normed sequence
space over the real or complex field. Let $q\in[1,\infty)$. We will
say that \textit{$E$ has property $S_{q}(C)$ }for some positive
constant $C$ if, whenever $\left\{ u_{n}\right\} _{n\in\mathbb{N}}$
and $\left\{ v_{n}\right\} _{n\in\mathbb{N}}$ are two sequences in
$\ell^{q}$ whose nonincreasing rearrangements $\left\{ u_{n}^{*}\right\} _{n\in\mathbb{N}}$
and $\left\{ v_{n}^{*}\right\} _{n\in\mathbb{N}}$ satisfy 
\begin{equation}
\sum_{n=1}^{N}(u_{n}^{*})^{q}\le\sum_{n=1}^{N}(v_{n}^{*})^{q}\mbox{ for all }N\in\mathbb{N}\label{eq:HL-one}
\end{equation}
 and 
\begin{equation}
\sum_{n=1}^{\infty}(u_{n}^{*})^{q}=\sum_{n=1}^{\infty}(v_{n}^{*})^{q}\label{eq:HL-two}
\end{equation}
and also $\left\{ u_{n}\right\} _{n\in\mathbb{N}}\in E$, then it
follows that $\left\{ v_{n}\right\} _{n\in\mathbb{N}}\in E$ and $\left\Vert \left\{ v_{n}\right\} _{n\in\mathbb{N}}\right\Vert _{E}\le C\left\Vert \left\{ u_{n}\right\} _{n\in\mathbb{N}}\right\Vert _{E}$.
\end{defn}

We are very indebted to Galina Levitina, Fedor Sukochev and Dmitriy
Zanin who motivated our research by having the remarkable insight
to realize that a property like $S_{q}(C)$ might characterize interpolation
spaces between $\ell^{q}$ and $\ell^{p}$ for some $p\in(0,q)$.
They express this explicitly in their Conjecture 1.5 on p.~3 of \cite{LeSuZa2017B}
which is formulated for the case $q=2$ and refers to quasi-normed
as well as normed sequence spaces. (See Remark \ref{rem:DetailedComparison}
for more details, comparing their conjecture with our results here.)

Their consideration of these matters arises from their research in
\cite{LeSuZa2017,LeSuZa2017B} into topics in operator theory which
also have some connections with mathematical physics. They give indications
of the plausibility of their conjecture in their Proposition 2.7 which
appears on p.~7 of \cite{LeSuZa2017} and on p.~8 of \cite{LeSuZa2017B}.
This proposition shows that, for suitable choices of the constant
$C$ and for all $p\in(0,2)$ the normed or quasi-normed spaces $L^{p}$
and $L^{p,\infty}$ (i.e., ``Weak $L^{p}$'') on $\left(0,\infty\right)$
both have a property for function spaces which is the exact analogue
of $S_{2}(C)$. A simple modification of their proof should show that
$\ell^{p}$ and $\ell^{p,\infty}$ have property $S_{2}(C)$ itself
for $p\in(0,2)$. Here, if $p<1$, we of course have to consider a
slight modification of Definition \ref{def:SqC} where $E$ may be
merely quasi-normed.

We have confined ourselves to considering only normed spaces in this
paper, but hope, in future research, to extend at least some of our
results to the more general setting of quasi-normed spaces. 

At a first quick glance it may seem very doubtful that property $S_{q}(C)$
could be relevant for describing interpolation spaces for the couple
$\left(\ell^{1},\ell^{q}\right)$. The sums appearing in (\ref{eq:HL-one})
are equivalent to the $K$-functionals (for $t=N$) of the sequences
$\left\{ u_{n}\right\} _{n\in\mathbb{N}}$ and $\left\{ v_{n}\right\} _{n\in\mathbb{N}}$
with respect to a \textit{different} couple, namely $\left(\ell^{q},\ell^{\infty}\right)$.
So $S_{q}(C)$ seems to be the condition which is known (essentially
by\cite{LorentzGShimogakiT1971}) to be equivalent to the space $E$
being an interpolation space with respect to $\left(\ell^{q},\ell^{\infty}\right)$,
except that $u_{n}^{*}$ and $v_{n}^{*}$ have been interchanged in
(\ref{eq:HL-one}). However these doubts vanish when one realizes
that (\ref{eq:HL-one}) can be re-expressed in terms of sums from
$N$ to $\infty$. (Analogously, in the proof of Proposition 2.7 in
\cite{LeSuZa2017,LeSuZa2017B} a condition in terms of integrals on
the interval $\left(0,t\right)$ is replaced by a condition in terms
of integrals on the interval $\left(t,\infty\right)$.) 

We can now begin the systematic presentation of our results by recalling
some terminology and known results about sequence spaces and from
interpolation of normed linear spaces.

Each definition and each result presented in this paper must be understood,
unless explicitly stated otherwise, to apply simultaneously to linear
spaces over either of the two base fields of real and complex numbers. 

All sequence spaces appearing here must be understood to be as indexed
using the natural numbers. Similarly, all $L^{p}$ spaces and other
spaces of (equivalence classes of) real or complex valued measurable
functions are to be understood, if not specified otherwise, to be
defined using $[0,\infty)$ equipped with Lebesgue measure as the
underlying measure space. 

All linear spaces considered here will be normed. However, as already
implied above, it seems reasonable to conjecture that it might be
possible to extend at least some parts of our results to the more
general realm of quasinormed spaces. 

We shall assume familiarity with the definitions and basic properties
of two standard objects, namely:

(i) the sequence $\left\{ x_{n}^{*}\right\} _{n\in\mathbb{N}}$ which
is the nonincreasing rearrangement of an arbitrary bounded sequence
$\left\{ x_{n}\right\} _{n\in\mathbb{N}}$ of real or complex numbers,

and also 

(ii) the function $f^{*}:(0,\infty)\to[0,\infty)$ which is the nonincreasing
rearrangement of an arbitrary element $f$ in the space $L^{1}(\mu)+L^{\infty}(\mu)$
defined on some arbitrary underlying measure space $\left(\Omega,\Sigma,\mu\right)$. 

In particular we will use the following very obvious connection between
(i) and (ii):
\begin{fact}
\label{fact:FstarXstar}If $\left\{ x_{n}\right\} _{n\in\mathbb{N}}$
is a real or complex sequence in $c_{0}$ and if $f$ is the function
on $[0,\infty)$ defined by $f=\sum_{n=1}^{\infty}x_{n}\chi_{[n-1,n)}$
then the function $f^{*}$ is given by the formula $f^{*}=\sum_{n=1}^{\infty}x_{n}^{*}\chi_{[n-1,n)}$.
\end{fact}
We shall also assume some familiarity with basic notions and results
in the theory of interpolation spaces. But it seems appropriate to
explicitly recall some of these notions and results, together with
the notation and terminology which we shall use for them here, which
will sometimes be slightly ``non-standard''. Some comments about
the classical interpolation theorems of Marcel Riesz and Riesz-Thorin
are made in an appendix (namely, Section \ref{sec:MarcelOlof}).
\begin{defn}
\label{def:ScriptL}For any (real or complex) Banach couples $\left(A_{0},A_{1}\right)$
and $\left(B_{0},B_{1}\right)$ we shall, as is often done in interpolation
theory, use the notation $T:\left(A_{0},A_{1}\right)\to\left(B_{0},B_{1}\right)$
to mean that $T$ is a linear map of $A_{0}+A_{1}$ into $B_{0}+B_{1}$
which satisfies $T\left(A_{j}\right)\subset B_{j}$ for $j=0,1$ and
also for which the quantity
\begin{align*}
\left\Vert T\right\Vert _{\left(A_{0},A_{1}\right)\to\left(B_{0},B_{1}\right)} & :=\max_{j=0,1}\left\{ \sup\left\Vert T\right\Vert _{A_{j}\to B_{j}}\right\} \\
 & :=\max_{j=0,1}\left\{ \sup\left\Vert Ta\right\Vert _{B_{j}}:a\in A_{j},\,\left\Vert a\right\Vert _{A_{j}}\le1\right\} 
\end{align*}
is finite. For each $r>0$ it will sometimes be convenient to let
$\mathcal{L}_{r}\left(\left(A_{0},A_{1}\right)\to\left(B_{0},B_{1}\right)\right)$
denote the set of all linear maps $T:\left(A_{0},A_{1}\right)\to\left(B_{0},B_{1}\right)$
which satisfy $\left\Vert T\right\Vert _{\left(A_{0},A_{1}\right)\to\left(B_{0},B_{1}\right)}\le r$.
In the case where $\left(B_{0},B_{1}\right)=\left(A_{0},A_{1}\right)$
we can use the abbreviated notation $\mathcal{L}_{r}\left(\left(A_{0},A_{1}\right)\right)$
for this set.
\end{defn}

\begin{defn}
\label{def:C-interp}Let $\left(A_{0},A_{1}\right)$ be a Banach couple,
let $C$ be a positive constant, and let $A\ne\left\{ 0\right\} $
be a normed space contained in $A_{0}+A_{1}$. We shall say that $A$
is a \textit{normed interpolation space with respect to the couple}
$\left(A_{0},A_{1}\right)$ if every $T\in\mathcal{L}_{1}\left(\left(A_{0},A_{1}\right)\right)$
has the property that $T(A)\subset A$. 

If, furthermore, 
\[
\left\Vert T\right\Vert _{A\to A}:=\sup\left\{ \left\Vert Ta\right\Vert _{A}:a\in A,\,\left\Vert a\right\Vert _{A}\le1\right\} \le C
\]
for each such $T$ then we shall say that $A$ is a \textit{normed
$C$-interpolation space with respect to the couple} $\left(A_{0},A_{1}\right)$.
In the case where $C=1$ we say that $A$ is an \textit{exact} normed
interpolation space. 

We can use the notation $Int\left(\left(A_{0},A_{1}\right)\right)$
for the collection of all normed interpolation spaces with respect
to $\left(A_{0},A_{1}\right)$ and $Int_{C}\left(\left(A_{0},A_{1}\right)\right)$
for the collection of all normed $C$-interpolation spaces with respect
to $\left(A_{0},A_{1}\right)$..\end{defn}
\begin{rem}
One can easily verify that any space $A$ which is a normed interpolation
space with respect to $\left(A_{0},A_{1}\right)$ must contain $A_{0}\cap A_{1}$.
We have made a point of inserting the word ``normed'' in the preceding
definition to emphasize the fact that here our terminology is slightly
different from that which is often used in the literature. The difference
is that here we are not requiring $A$ to be complete, i.e.,~to be a Banach
space, nor requiring that the stated inclusion $A\subset A_{0}+A_{1}$
nor the inclusion $A_{0}\cap A_{1}\subset A$ be continuous. We permit
ourselves to make this slight deviation from standard definitions
and terminology since it is just as easy to obtain our results for
this larger class of interpolation spaces as it is for their subset
of those which have the additional properties just mentioned. (Some
further non-essential comments about this appear in 
Subsection \ref{subsec:DiffDefnInterpSp}.)\end{rem}
\begin{defn}
\label{def:Kspace}Let 
$\left(A_{0},A_{1}\right)$ be a Banach couple
and let $C$ be a positive constant. A normed space 
$A\ne\left\{ 0\right\}$
which is contained in $A_{0}+A_{1}$ is said to be a 
\textit{normed $K$ 
space with respect to }$\left(A_{0},A_{1}\right)$ if, whenever
$x$ and $y$ are elements of $A_{0}+A_{1}$ whose $K$-functionals
satisfy
\begin{equation}
K(t,y;A_{0},A_{1})\le K(t,x;A_{0},A_{1})\mbox{\,\ for all\,}t>0\label{eq:prink}
\end{equation}
and $x$ is an element of $A$, then it follows that $y\in A$. If,
furthermore, the preceding conditions on $x$ and $y$ also imply
that $\left\Vert y\right\Vert _{A}\le C\left\Vert x\right\Vert _{A}$
then we say that $A$ is a\textit{ normed $C-K$ space with respect
to }$\left(A_{0},A_{1}\right)$. In the case where $C=1$ we say that
$A$ is an \textit{exact} normed $K$ space. 
\end{defn}
As in Definition \ref{def:C-interp}, we have inserted the word ``normed''
also in this definition to emphasize that the space $A$, is merely
assumed to be normed, and not necessarily complete. Nor are the inclusions
$A_{0}\cap A_{1}\subset A\subset A_{0}+A_{1}$ required to be continuous.
Here too we relegate further comments about this to Subsection 
\ref{subsec:DiffDefnInterpSp}.
\begin{rem}
\label{rem:BiggerThanOne-1}In the context of Definition \ref{def:C-interp},
we note that obviously the identity operator on $A_{0}+A_{1}$ is
in $\mathcal{L}_{1}\left(\left(A_{0},A_{1}\right)\right)$. In the
context of Definition \ref{def:Kspace}, we note that obviously (\ref{eq:prink})
holds when $x=y$. From these two observations we see that any constant
$C$ which appears in Definition \ref{def:C-interp} or in Definition
\ref{def:Kspace} must satisfy $C\ge1$. 
\end{rem}
Standard properties of the $K$-functional immediately imply the following
frequently used fact (though it is often stated only in the context
of complete interpolation spaces and complete $K$-spaces).
\begin{fact}
\label{fac:KspaceIsInterp}For every Banach couple $\left(A_{0},A_{1}\right)$
every $K$ space with respect to $\left(A_{0},A_{1}\right)$ is also
an interpolation space with respect to $\left(A_{0},A_{1}\right)$.
Furthermore, for every $C\ge1$, every space $A$ which is a normed
$C-K$ space with respect to $\left(A_{0},A_{1}\right)$ is also a
normed $C$-interpolation space with respect to $\left(A_{0},A_{1}\right)$.
\end{fact}
It is also very well known that, for some special couples, among them,
the couple $(L^{1}(\mu),L^{\infty}(\mu))$ for any%
\footnote{For the extension of this and other results to the case where the
measure space is not $\sigma$-finite see \cite[pp.~232-233]{CwikelM1976}.%
}underlying measure space $\left(\Omega,\Sigma,\mu\right)$ (see \cite{CalderonA1966},
cf.~also \cite{MityaginB1965}), a converse of Fact \ref{fac:KspaceIsInterp}
also holds. Essentially as in Part (iii) of Theorem 3 on p.~208 of
\cite{CalderonA1966}, we have 
\begin{thm}
\label{thm:CalderonExplicit}, Every space $A$ contained in $L^{1}(\mu)+L^{\infty}(\mu)$
is a normed exact interpolation space with respect to $\left(L^{1}(\mu),L^{\infty}(\mu)\right)$
if and only if it is a normed exact K space with respect to $\left(L^{1}(\mu),L^{\infty}(\mu)\right)$.
\end{thm}
(Again completeness of $A$ and continuity of the inclusions $A_{0}\cap A_{1}\subset A\subset A_{0}+A_{1}$
are each irrelevant for showing this.) Since there is an exact ``concrete''
formula for the $K$-functional of $(L^{1}(\mu),L^{\infty}(\mu))$,
Theorem \ref{thm:CalderonExplicit} is very easily seen (cf.~also
Lemma 1 on p.~232 of \cite{SparrG1978}) to be an immediate consequence%
\footnote{In fact a simple argument using ``orbit spaces'' gives almost the
reverse implication, namely that Theorem \ref{thm:CalderonExplicit}
also implies a slightly weaker form of Theorem \ref{thm:Calderon}
in which, for any given positive $\varepsilon$, the norm $\left\Vert T\right\Vert _{\mathcal{L}\left((L^{1}(\mu),L^{\infty}(\mu))\right)}$
can only be asserted to be less than $1+\varepsilon$. This is a special
case of a result which can be formulated for arbitrary Banach couples.
Cf., e.g., Section 4.1 on p.~208 of \cite{NilssonP1983} or, qualitatively,
Lemma 4.4.7 on p.~580 of \cite{BrudnyiYKrugljakN1991}. %
} of 
\begin{thm}
\label{thm:Calderon}\cite[Theorem 1, p.~278]{CalderonA1966} Let
$f$ and $g$ be elements of $L^{1}(\mu)+L^{\infty}(\mu)$ whose nonincreasing
rearrangements $f^{*}$ and $g^{*}$ satisfy 
\[
\int_{0}^{t}g^{*}(s)ds\le\int_{0}^{t}f^{*}(s)ds\mbox{\,\ for every\,}t>0.
\]
Then there exists a linear map $T\in\mathcal{L}_{1}\left((L^{1}(\mu),L^{\infty}(\mu))\right)$
for which $Tf=g$.
\end{thm}
The above mentioned (well known!) formula for the $K$-functional
for this couple, for each $f\in L^{1}(\mu)+L^{\infty}(\mu)$, is of
course
\begin{equation}
K(t,f;L^{1}(\mu),L^{\infty}(\mu))=\int_{0}^{t}f^{*}(s)ds\mbox{\,\ for all\,}t>0\label{eq:KfunLoneLinf}
\end{equation}
and is due, independently, to Evelio Oklander \cite{OklanderE1964}
and to Jaak Peetre \cite{PeetreJ1968brasilia}.

For some couples, among them couples of $L^{p}$ spaces (see the paper
\cite{SparrG1978} of Gunnar Sparr, as well as the earlier special
cases of his results cited there), a ``qualitative'' converse of
Fact \ref{fac:KspaceIsInterp} holds, i.e., every $C$-interpolation
space for the couple is also a $rC-K$ space for some constant $r$
depending only on the couple. We shall not need to use Sparr's result
here. However, a variant of the approach appearing in its alternative
proofs in \cite{ArazyJCwikelM1984} and \cite{CwikelM1976} and also
one of the lemmata in \cite{ArazyJCwikelM1984} will play roles for
us.
\begin{rem*}
In general, the converse of Fact \ref{fac:KspaceIsInterp} does not
hold, not even qualitatively. I.e., for many couples $\left(A_{0},A_{1}\right)$
there exist spaces $A$ which are normed (and in fact complete)$C_{1}$-interpolation
spaces for some $C_{1}>0$ but which are not normed $C_{2}-K$ spaces
for any $C_{2}>0$. See e.g., \cite{CwikelM-Encyclopedia}. 
\end{rem*}
We will need to use a well known result which is a simple special
case of an argument which appears as Lemma 2 on pp.~277-278 of \cite{CalderonA1966}
where it is one of the steps in the proof of the result which we recalled
and formulated above as Theorem \ref{thm:Calderon}. It can also be
readily seen to be a special case of Theorem \ref{thm:Calderon}.
Let us state it explicitly here and, for convenience, briefly provide
a simple self-contained proof which bypasses the more elaborate reasoning
in \cite{CalderonA1966}.
\begin{lem}
\label{lem:EasyCalderon} Let $x=\left\{ x_{n}\right\} _{n\in\mathbb{N}}$
be a sequence of real or complex numbers which tends to $0$ and let
$x^{*}=\left\{ x_{n}^{*}\right\} _{n\in\mathbb{N}}$ be its nonincreasing
rearrangement. Then there exist linear maps $W$ and $Y$ in $\mathcal{L}_{1}\left(\left(\ell^{1},\ell^{\infty}\right)\right)$
for which $Wx=x^{*}$ and $Yx^{*}=x$.
\end{lem}
\noindent \textit{Outline of the proof.} Let us first explicitly
state a fact which enables the proof of this lemma to be much simpler
than that of \cite[Lemma 2, p.~277]{CalderonA1966}:

\textit{For each $\alpha>0$ the two sets $\left\{ n\in\mathbb{N}:\left|x_{n}\right|=\alpha\right\} $
and $\left\{ n\in\mathbb{N}:x_{n}^{*}=\alpha\right\} $ both have
the same finite cardinality (including the possibility that they are
both empty). }

Let $\Xi=\left\{ n\in\mathbb{N}:x_{n}\ne0\right\} $. If $\Xi$ is
finite, or is all of $\mathbb{N}$, then it is easy to see that we
can choose $W$ to be of the form 
\begin{equation}
W\left(\left\{ u_{n}\right\} _{n\in\mathbb{N}}\right)=\left\{ \theta_{n}u_{\sigma(n)}\right\} _{n\in\mathbb{N}}\mbox{ for each\,}\left\{ u_{n}\right\} _{n\in\mathbb{N}}\in\ell^{\infty},\label{eq:thfo}
\end{equation}
for some suitable sequence $\left\{ \theta_{n}\right\} _{n\in\mathbb{N}}$
of real or complex numbers satisfying $\left|\theta_{n}\right|=1$
for all $\mathbb{N}$ and for some suitable one to one map $\sigma$
of $\mathbb{N}$ onto itself. Otherwise, if $\Xi$ is an infinite
but proper subset of $\mathbb{N}$ we can still choose $W$ to be
given by a formula of the form (\ref{eq:thfo}) for a suitable sequence
$\left\{ \theta_{n}\right\} _{n\in\mathbb{N}}$ as above, but this
time $\sigma$ has to be a suitable one to one map of $\mathbb{N}$
onto $\Xi$. 

It is also easy to see that the map $Y$ can be given by a formula
similar to (\ref{eq:thfo}) in the cases where $\Xi$ is finite or
is all of $\mathbb{N}$. In these cases $\sigma$ is again a suitable
one to one map of $\mathbb{N}$ onto itself and $\left\{ \theta_{n}\right\} _{n\in\mathbb{N}}$
is a suitable sequence satisfying $\left|\theta_{n}\right|=1$ for
all $n$. This time, in the case where $\Xi$ is an infinite but proper
subset of $\mathbb{N}$, we choose $Y\left(\left\{ u_{n}\right\} _{n\in\mathbb{N}}\right)$
for each given $\left\{ u_{n}\right\} _{n\in\mathbb{N}}\in\ell^{\infty}$,
to be the sequence $\left\{ y_{n}\right\} _{n\in\mathbb{N}}$ defined
as follows: For each $n\in\mathbb{N}\setminus\Xi$ we set $y_{n}=0$
and, for each $n\in\Xi$ we set $y_{n}=\theta_{n}u_{\sigma(n)}$ for
suitable $\theta_{n}$ satisfying $\left|\theta_{n}\right|=1$ and
where $\sigma$ is a suitable one to one map of $\Xi$ onto $\mathbb{\mathbb{N}}$.
$\qed$

\medskip{}

\begin{rem}
\label{rem:NormIsOneForEveryR}Obviously the simple forms of the operators
$W$ and $Y$ in the above proof show, without our even needing to
apply the Riesz-Thorin or Marcel Riesz interpolation theorems, that
these operators also satisfy $\left\Vert W\right\Vert _{\ell^{r}\to\ell^{r}}=1$
and $\left\Vert Y\right\Vert _{\ell^{r}\to\ell^{r}}=1$ for all $r\in[1,\infty]$. 
\end{rem}
There will be several steps in some of our proofs in which it will
be convenient to use the two simple special linear operators $P$
and $Q$ which we introduce in the following statement:
\begin{fact}
\label{fact:PQ}Let the two linear maps $P:\ell^{\infty}\to L{}^{1}+L^{\infty}$
and $Q:L^{1}+L^{\infty}\to\ell^{\infty}$ be defined by 
\[
P\left(\left\{ h_{n}\right\} _{n\in\mathbb{N}}\right)=\sum_{n=1}^{\infty}h_{n}\chi_{[n-1,n)}
\]
for each $\left\{ h_{n}\right\} _{n\in\mathbb{N}}\in\ell^{\infty}$
and 
\[
Qh=\left\{ \int_{n-1}^{n}h(s)ds\right\} _{n\in\mathbb{N}}
\]
for each $h\in L^{1}+L^{\infty}$. 

Then, for each $r\in[1,\infty]$, $\left\Vert P\left(\left\{ h_{n}\right\} _{n\in\mathbb{N}}\right)\right\Vert _{L^{r}}=\left\Vert \left\{ h_{n}\right\} _{n\in\mathbb{N}}\right\Vert _{\ell^{r}}$
for each $\left\{ h_{n}\right\} _{n\in\mathbb{N}}\in\ell^{r}$ and
$\left\Vert Qh\right\Vert _{\ell^{r}}\le\left\Vert h\right\Vert _{L^{r}}$
for each $h\in L^{r}$, 

So of course $P\in\mathcal{L}_{1}\left((\ell^{r_{0}},\ell^{r_{1}})\to(L^{r_{0}},L^{r_{1}})\right)$
and $Q\in\mathcal{L}_{1}\left((L^{r_{0}},L^{r_{1}})\to(\ell^{r_{0}},\ell^{r_{1}})\right)$
for any $r_{0}$ and $r_{1}$ in $\left[1,\infty\right]$.
\end{fact}
\bigskip{}

\section{\label{sec:MainResults}The main result}

Let us first briefly mention some previous results about the couple
$\left(\ell^{1},\ell^{q}\right)$ and couples related to it. Apparently
G.~G.~Lorentz and Tetsuya Shimogaki were the first to obtain a characterization
of all the Banach function spaces with the Fatou property%
\footnote{Lorentz and Shimogaki express the fact that a given function space
has the Fatou property by using the alternative terminology that its
``norm is semicontinuous''. %
}which are interpolation spaces with respect to the couple $\left(L^{1},L^{q}\right)$
for $q\in(1,\infty)$. They expressed this result as their Theorem
3 on p.~216 of \cite{LorentzGShimogakiT1971} which is the second
of their two main results in that paper. From that result one can
deduce a related characterization of interpolation spaces with respect
to $\left(\ell^{1},\ell^{q}\right)$.

Above, in Section \ref{sec:Intro}, we already briefly mentioned the
results which Gunnar Sparr \cite{SparrG1978} obtained some seven
years after the work of Lorentz and Shimogaki. They include, as a
special case, a characterization of all interpolation spaces with
respect to the couple $\left(L^{1},L^{q}\right)$. Sparr's characterization,
as already indicated above, is explicitly in terms of the $K$-functional
for that couple. So it is quite different from Theorem 3 of \cite{LorentzGShimogakiT1971},
and instead is analogous%
\footnote{The authors of \cite{CalderonA1966} and \cite{LorentzGShimogakiT1971}
were apparently unaware at the time of writing those papers that the
functionals that they were using were equivalent to the $K$-functionals
of the couples that they were studying. So this analogy became apparent
only later.%
}to the result of Calder\'on for $\left(L^{1},L^{\infty}\right)$
in \cite{CalderonA1966} mentioned above, and in fact to the other
main result in \cite{LorentzGShimogakiT1971} which applies to the
couple $\left(L^{p},L^{\infty}\right)$). Sparr's characterization
also applies to spaces which do not necessarily have the Fatou property.
Furthermore, it is valid for all choices of the underlying measure
space and therefore it requires no modification to make it immediately
applicable to the couple $\left(\ell^{1},\ell^{q}\right)$ which we
are considering here.

The main notion in the formulation of our new rather different
way of characterizing interpolation spaces for the couple $\left(\ell^{1},\ell^{q}\right)$
is the property $S_{q}(C)$ which we already introduced above in Definition
\ref{def:SqC}. But we also need another notion, namely the following
weakened version of the standard Fatou property. 
\begin{defn}
\label{def:WeakFatou}For each sequence $\left\{ x_{n}\right\} _{n\in\mathbb{N}}$
of real or complex numbers and for each $N\in\mathbb{N}$, let $\left\{ x_{n}^{(N)}\right\} _{n\in\mathbb{N}}$
be the truncated sequence defined by $x_{n}^{(N)}=x_{n}$ for all
$n\in\left\{ 1,2,...,N\right\} $ and $x_{n}^{(N)}=0$ for all $n>N$.
Let $R$ be a positive constant. Let $E$ be a normed sequence space.
Suppose that whenever $\left\{ x_{n}\right\} _{n\in\mathbb{N}}$ is
a sequence of non-negative numbers for which $\left\{ x_{n}^{(N)}\right\} _{n\in\mathbb{N}}\in E$
for all $N\in\mathbb{N}$ and $\sup_{N\in\mathbb{N}}\left\Vert \left\{ x_{n}^{(N)}\right\} _{n\in\mathbb{N}}\right\Vert _{E}<\infty$,
then it follows that $\left\{ x_{n}\right\} _{n\in\mathbb{N}}\in E$
and 
\[
\left\Vert \left\{ x_{n}\right\} _{n\in\mathbb{N}}\right\Vert _{E}\le R\sup_{N\in\mathbb{N}}\left\Vert \left\{ x_{n}^{(N)}\right\} _{n\in\mathbb{N}}\right\Vert _{E}.
\]
Then we shall say that $E$ has property $WFP(R).$
\end{defn}
The initials WFP can be thought of as standing for ``Weak Fatou Property'',
although this latter terminology is sometimes given a slightly different
meaning in the literature. There are other equivalent definitions
of this property, especially for when $E$ is known to also have other
properties. We have chosen what seems to be the minimal requirement
which stipulates the property that we require for our purposes below
in Theorem \ref{thm:Enough}.

\begin{rem}
\label{rem:BiggerThanOne}By considering sequences with finite support
we see that any constant $R$ appearing in Definition \ref{def:WeakFatou}
must satisfy $R\ge1$. By considering what happens when the sequences
$\left\{ u_{n}\right\} _{n\in\mathbb{N}}$ and $\left\{ v_{n}\right\} _{n\in\mathbb{N}}$
in Definition \ref{def:SqC} are the same sequence, we also see that
any constant $C$ which appears in Definition \ref{def:SqC} must
satisfy $C\ge1$. 
\end{rem}
The main result is stated as:
\begin{thm}
\label{thm:ReallyMain}Let $q\in[1,\infty)$ and let $E$ be a normed
sequence space which is contained in $\ell^{q}$ and which has the
property $WFP(R)$ for some constant $R\ge1$. Then $E$ is a $C$-interpolation
space with respect to $\left(\ell^{1},\ell^{q}\right)$ for some constant
$C\ge1$ if and only if it has property $S_{q}(C')$ for some possibly
different constant $C'\ge1$.

More precisely: 

$\mathrm{(i)}$ If $E$ is a normed $C$-interpolation space with
respect to $\left(\ell^{1},\ell^{q}\right)$ for some constant $C\ge1$
then it has property $S_{q}(C)$, 

and

$\mathrm{(ii)}$ If $E$ has property $S_{q}(C)$ for some constant
$C\ge1$ then it is a $C'$-interpolation space and also a $C'-K$
space with respect to $\left(\ell^{1},\ell^{q}\right)$ for a certain
constant $C'\ge1$ which depends only on $R$, $C$ and $q$. 
\end{thm}
\noindent \textit{Proof.} This theorem follows from the combined
results of Theorems \ref{thm:Easy}, \ref{thm:Main} and \ref{thm:Enough}
below. $\qed$
\begin{rem}
\label{rem:VeryImportant} This theorem and also the related Theorems
\ref{thm:Easy}, \ref{thm:Main} and \ref{thm:Enough}, all apply,
in particular, to the special case where the normed sequence space
$E$ is complete and also satisfies the continuous inclusions $\ell^{1}\subset E\subset\ell^{q}$.
Then the additional properties are unnecessary for the proofs of these
theorems (though the continuity of the above mentioned inclusions
may be a consequence of other hypotheses on $E$). So our results
also cover the case of intermediate spaces and $C$-interpolation
spaces when they are defined in the usual way.
\end{rem}
In fact, as we shall now specify, Theorems \ref{thm:Easy}, \ref{thm:Main}
and \ref{thm:Enough} also provide more information than is expressed
by Theorem \ref{thm:ReallyMain}. We cannot dispense with the requirement
that the normed sequence space $E$ is contained in $\ell^{q}$, but
some parts of Theorem \ref{thm:ReallyMain} hold when the requirement
that $E$ has property $WFP(R)$ is dispensed with, or replaced by
a less stringent requirement. 

$\bullet$ The implication in part (i) of Theorem \ref{thm:ReallyMain}
follows from Theorem \ref{thm:Easy} and that theorem shows that it
holds also when $E$ does not have the property $WFP(R)$ for any
$R$. 

$\bullet$ The combined results of Theorems \ref{thm:Enough} and
\ref{thm:Main} establish the implication in part (ii) of Theorem
\ref{thm:ReallyMain}. In fact these theorems show that it also holds
when, in addition to our standing general requirement that $E$ is
contained in $\ell^{q}$ and that (in part (ii)) $E$ has property
$S_{q}(C)$, instead of requiring $E$ to have the property $WFP(R)$
for some $R$, we require it to have a different property, which in
our context is weaker than $WFP(R)$, That different property is that
$E$ is a $C_{2}$-interpolation space with respect to the couple
$\left(\ell^{1},\ell^{\infty}\right)$ for some constant $C_{2}\ge1$.
More precisely, it is Theorem \ref{thm:Enough} which shows that this
latter property really is weaker in our context, i.e., that $WFP(R)$
\textit{together with containment in $\ell^{q}$ and property $S_{q}(C)$
}imply that $E\in Int_{C_{2}}\left((\ell^{1},\ell^{\infty})\right)$
for $C_{2}=CR$. It is Theorem \ref{thm:Main} which shows that, even
if the normed sequence space $E$ happens not to have property $WFP(R)$
for any $R$, but does satisfy $E\in Int_{C_{2}}\left((\ell^{1},\ell^{\infty})\right)$,
as well as the property $S_{q}(C)$ which is imposed in part (ii)
of Theorem \ref{thm:ReallyMain}, this suffices to ensure that $E$
is a $C'-K$ space and therefore also a $C'$-interpolation space
with respect to the couple $(\ell^{1},\ell^{q})$ for some constant
$C'$ depending only on $q$, $C$ and $C_{2}$. 
\begin{rem}
In the case $q=1$ the pair of inequalities (\ref{eq:HL-one}) and
(\ref{eq:HL-two}), play a role in results going back to Hardy, Littlewood
and P\'olya. For further historical remarks and other applications
related to (\ref{eq:HL-one}) and (\ref{eq:HL-two}), see, e.g.,~Theorem
HLP on pp.~233\textbf{-}234 of \cite{SparrG1978} and the subsection
entitled ``Majorization'' on pp.~6-9 of \cite{LeSuZa2017}. We
shall see later that these two inequalities imply other inequalities
(cf.~(\ref{eq:rtq}) and (\ref{eq:aqe}) in the proof below of Theorem
\ref{thm:Easy}) between sums or integrals which are terms in Tord
Holmstedt's equivalent formula (to which we will refer in more detail
below in Section \ref{sec:Main}, see (\ref{eq:Holmstedt})) for the
$K$-functional for $\left(\ell^{1},\ell^{q}\right)$ or for $\left(L^{1},L^{q}\right)$.
But they are not the whole of that formula, and so apparently none
of our results are trivial consequences of Sparr's result.
\end{rem}

\begin{rem}
\label{rem:DetailedComparison}To conclude this section we offer an
explicit comparison of our results with Conjecture 1.5 of \cite[p.~3]{LeSuZa2017B}.
Expressed in our teminology, the statement of Conjecture 1.5 is: 

\textit{Let $E$ be a normed or quasi-normed symmetric sequence space.
Then the following conditions are equivalent:}

\textit{$\mathrm{(i)}$~$E$ has property $S_{2}(C)$ for some positive
constant $C$.}

\textit{$\mathrm{(ii)}$ There exists some $p\in(0,2)$ such that
$E$ is a interpolation space with respect to $\left(\ell^{p},\ell^{2}\right)$.}

We first note that condition (ii) obviously implies that $E\subset\ell^{2}$,
and that the authors of \cite{LeSuZa2017B} have clarified elsewhere
(see the paragraph just before Remark \ref{rem:ExoticExamples}) that
condition (i) also implies that $E\subset\ell^{2}$. 

Though we hope to extend them in future research, our results here
only apply to the case where $E$ is normed and the number $p$ in
condition (ii) is in the range $[1,2)$. Subject to these restrictions
we can prove that slight variants of conditions (i) and (ii) imply
each other:

If we assume that $E$ satisfies condition (ii) and, furthermore,
that it is not only an interpolation space with respect to the couple
$\left(\ell^{p},\ell^{2}\right)$ but also a $C$-interpolation space
for some $C$ with respect to that couple , then our Theorem \ref{thm:Easy}
(in the particular case where $q=2$) enables us to deduce that condition
(i) also holds. We also obtain that the constant $C$ in property
$S_{2}(C)$ of $E$ turns out to be in fact the same constant as in
our assumption about $E$. Note that if $E$ is complete then we obtain
``for free'' by \cite[Theorem 6.XI, p.~73]{AronszajnNGagliardoE1965}
that it is also a $C$-interpolation space for some $C$. 

As discussed at the beginning of the proof of Theorem \ref{thm:Easy},
among all values of $p$ in {[}1,2), the condition that $E$ is a
$C$-interpolation space with respect to $\left(\ell^{p},\ell^{2}\right)$
is weakest when $p=1$. We point out that our above mentioned variant
of the implication $\mathrm{(ii)\Rightarrow\mathrm{(i)}}$ holds when
we require our version of condition $\mathrm{(ii)}$ to hold merely
for this least demanding choice of $p$. 

Our version of the reverse implication $\mathrm{(i)}\Rightarrow\mathrm{(ii)}$
is as follows: Suppose that our normed symmetric sequence space $E$
satisfies condition (i) i.e., that it has property $S_{2}(C)$. Then,
subject to one further requirement, we can deduce that condition (ii)
holds, where the exponent $p$ in that condition can always be chosen
to equal $1$. (Note that for some choices of $E$, among them $\ell^{1}$,
no other choice of $p$ in $[1,2)$ is possible.) That one further
requirement can either be that $E$ has the Fatou property or a weaker
version of that property (see Definition \ref{def:WeakFatou}), or
that $E$ is a normed $C'$-interpolation space for some constant
$C'$ with respect to the couple $\left(\ell^{1},\ell^{\infty}\right)$.
Our conclusion in fact a little stronger than condition (ii) for $p=1$,
namely, we obtain that $E$ is a $C^{\prime\prime}$-interpolation
space with respect to $\left(\ell^{1},\ell^{2}\right)$ where the
constant $C^{\prime\prime}$ can be bounded by an explicit quantity
depending only on $C$, $C'$ and possibly also on a constant which
arisess in the definition of the weak Fatou property for $E$. We
obtain this result by applying our Theorem \ref{thm:Main} in the
case $q=2$ and also, if necessary, applying our Theorem \ref{thm:Enough}
for $q=2$.
\end{rem}

\section{Interpolation spaces between $\ell^{1}$ and $\ell^{q}$ have property
$S_{q}(C)$}
\begin{thm}
\label{thm:Easy}Suppose that the constants $p$, $q$ and $C$ satisfy
$1\le p<q<\infty$ and $C>0$. Let $E$ be a normed $C$-interpolation
space with respect to the couple $(\ell^{p},\ell^{q})$. Then $E$
has property $S_{q}(C)$.\end{thm}
\begin{rem}
\label{rem:CompareWithFSandCoauthors} Among our theorems this is
the one which is closest to the result of Proposition 2.7 on p.~7
of \cite{LeSuZa2017} and on p.~8 of \cite{LeSuZa2017B}, to which
we already referred above, and which gives a similar result for the
case where $q=2$ where the underlying measure space is $\left(0,\infty\right)$
with Lebesgue measure instead of $\mathbb{N}$ with counting measure.
More explicitly, Proposition 2.7 shows that for each $r\in(0,2)$
the particular Banach or quasi-Banach spaces $L^{r}(0,\infty)$ and
$L^{r,\infty}(0,\infty)$ both have a property analogous to $S_{2}(C)$
for some constant $C$ depending on $r$, and in fact $C=1$ for that
property for $L^{r}(0,\infty)$. 
\end{rem}
\noindent \textit{Proof of Theorem \ref{thm:Easy}.} By one or the
other of the classical theorems of Marcel Riesz and of Riesz-Thorin
(cf.~Section \ref{sec:MarcelOlof}), the collection of operators
$\mathcal{L}_{1}\left(\left(\ell^{1},\ell^{q}\right)\right)$ is contained
in the collection $\mathcal{L}_{1}\left(\left(\ell^{p},\ell^{q}\right)\right)$
for each $p\in(1,q)$. So, if $E$ is a normed $C$-interpolation
space with respect to $\left(\ell^{p},\ell^{q}\right)$ for some $p\in(1,q)$,
then it is also a normed $C$-interpolation space with respect to
$\left(\ell^{1},\ell^{q}\right)$. This enables us to assume in this
proof, without loss of generality, that $p=1$.

Suppose that the sequences $x=\left\{ x_{n}\right\} _{n\in\mathbb{N}}$
and $y=\left\{ y_{n}\right\} _{n\in\mathbb{N}}$ of non-negative numbers
are non-increasing and satisfy 
\begin{equation}
\sum_{n=1}^{N}x_{n}^{q}\le\sum_{n=1}^{N}y_{n}^{q}\mbox{ for all }N\in\mathbb{N}\label{eq:zHL-one}
\end{equation}
 and 
\begin{equation}
\sum_{n=1}^{\infty}x_{n}^{q}=\sum_{n=1}^{\infty}y_{n}^{q}<\infty.\label{eq:Equals}
\end{equation}
Then 
\[
\sum_{n=1}^{\infty}y_{n}^{q}-\sum_{n=1}^{N}y_{n}^{q}\le\sum_{n=1}^{\infty}x_{n}^{q}-\sum_{n=1}^{N}x_{n}^{q}\mbox{ for every }N\in\mathbb{N},
\]
and this together with (\ref{eq:Equals}) implies that 
\begin{equation}
\sum_{n=M}^{\infty}y_{n}^{q}\le\sum_{n=M}^{\infty}x_{n}^{q}<\infty\mbox{ for every }M\in\mathbb{N}.\label{eq:rtq}
\end{equation}
Consider the two non-increasing functions $f:[0,\infty)\to[0,\infty)$
and $g:[0,\infty)\to[0,\infty)$ defined by $f=\sum_{n=1}^{\infty}x_{n}\chi_{[n-1,n)}$
and $g=\sum_{n=1}^{\infty}y_{n}\chi_{[n-1,n)}$. Obviously (cf.~Fact
\ref{fact:FstarXstar}) we can reformulate the condition (\ref{eq:rtq})
by stating that 
\begin{equation}
\int_{t}^{\infty}g(s)^{q}ds\le\int_{t}^{\infty}f(s)^{q}ds\label{eq:aqe}
\end{equation}
for every non-negative integer $t$. Since the functions of $t$ on
both sides of (\ref{eq:aqe}) are affine on the interval $\left[n-1,n\right]$
for $n\in\mathbb{N}$, we deduce that in fact (\ref{eq:aqe}) holds
for every $t\ge0$. So this is precisely the inequality which appears
and is labelled as (5) in the statement of the lemma on pp.~257-258
of \cite{ArazyJCwikelM1984}. Since the functions $f$ and $g$ are
both non-increasing, non-negative and countably valued, that lemma
can be applied to ensure the existence of a linear operator $V:(L^{1},L^{q})\to(L^{1}.L^{q})$
with $\left\Vert V\right\Vert _{(L^{1},L^{q})\to(L^{1},L^{q})}\le1$
for which $Vf=g$. (Perhaps we can say more about $V$. See Remark
\ref{rem:ArazyCwikelConjecture} below.)

We will now use $V$ in a very obvious way to obtain a linear map
$V_{1}:(\ell^{1},\ell^{q})\to(\ell^{1},\ell^{q})$ which satisfies
\begin{equation}
\left\Vert V_{1}\right\Vert _{(\ell^{1},\ell^{q})\to(\ell^{1},\ell^{q})}\le1\mbox{ ,}\label{eq:BothV}
\end{equation}
and, furthermore, $V_{1}x=y$. We simply take $V_{1}$ to be the composed
operator $V_{1}=QVP$, where, as in Fact \ref{fact:PQ}, $P:\ell^{\infty}\to L^{1}+L^{\infty}$
is defined by $P\left(\left\{ h_{n}\right\} _{n\in\mathbb{N}}\right)=\sum_{n=1}^{\infty}h_{n}\chi_{[n-1,n)}$
for every $\left\{ h_{n}\right\} _{n\in\mathbb{N}}$ in $\ell^{\infty}$
and $Q:L^{1}+L^{\infty}\to\ell^{\infty}$ is defined by $Qh=\left\{ \int_{n-1}^{n}h(s)ds\right\} _{n\in\mathbb{N}}$
for every $h\in L^{1}+L^{\infty}$. As already noted in Fact \ref{fact:PQ},
obviously $\left\Vert P\right\Vert _{\ell^{r}\to L^{r}}=1$ and $\left\Vert Q\right\Vert _{L^{r}\to\ell^{r}}\le1$
for all $r\in[1,\infty]$. So we obtain that (\ref{eq:BothV}) holds.
It is also obvious that $Px=f$ and $Qg=y$ and, consequently, $V_{1}x=y$.

. Our hypotheses imply that $E\subset\ell^{q}$ and also enable us
to renorm $E$ in a standard way to make it an exact interpolation
space, i.e., an element of $Int_{1}\left((\ell^{1},\ell^{q})\right)$.
We simply set
\[
\left\Vert x\right\Vert _{\widetilde{E}}:=\sup_{T\in\mathcal{L}_{1}\left((\ell^{1},\ell^{q})\right)}\left\Vert Tx\right\Vert _{E}\mbox{ for each\,}x\in E.
\]
Our hypotheses on $E$ obviously ensure that the quantity $\left\Vert \cdot\right\Vert _{\widetilde{E}}$
is a norm on $E$ which satisfies 
\[
\left\Vert x\right\Vert _{E}\le\left\Vert x\right\Vert _{\widetilde{E}}\le C\left\Vert x\right\Vert _{E}\mbox{ for all\,}x\in E
\]
and also 
\begin{equation}
Tx\in E\mbox{ with\,}\left\Vert Tx\right\Vert _{\widetilde{E}}\le\left\Vert x\right\Vert _{\widetilde{E}}\mbox{ for all\,}x\in E\mbox{ and all\,}T\in\mathcal{L}_{1}\left((\ell^{1},\ell^{q})\right).\label{eq:ryz}
\end{equation}
 Let $\left\{ u_{n}\right\} _{n\in\mathbb{N}}$ and $\left\{ v_{n}\right\} _{n\in\mathbb{N}}$
be sequences in $\ell^{q}$ whose nonincreasing rearrangements satisfy
(\ref{eq:HL-one}) and (\ref{eq:HL-two}) and suppose, furthermore,
that $\left\{ u_{n}\right\} _{n\in\mathbb{N}}\in E$. Since $\lim_{n\to\infty}u_{n}=0$
we can invoke Lemma \ref{lem:EasyCalderon} and Remark \ref{rem:NormIsOneForEveryR}
to obtain an operator $W$ which is a norm $1$ map of $\ell^{r}$
into itself for each $r\in[1,\infty]$ and for which $W\left(\left\{ u_{n}\right\} _{n\in\mathbb{N}}\right)$
is the nonincreasing rearrangement $\left\{ u_{n}^{*}\right\} _{n\in\mathbb{N}}$
of $\left\{ u_{n}\right\} _{n\in\mathbb{N}}$. 

Consequently we can apply (\ref{eq:ryz}) to obtain that $\left\{ u_{n}^{*}\right\} _{n\in\mathbb{N}}\in E$
and 
\[
\left\Vert \left\{ u_{n}^{*}\right\} _{n\in\mathbb{N}}\right\Vert _{\widetilde{E}}\le\left\Vert \left\{ u_{n}\right\} _{n\in\mathbb{N}}\right\Vert _{\widetilde{E}}.
\]

We can now apply the reasoning of the first part of this proof to
the two particular sequences $\left\{ x_{n}\right\} _{n\in\mathbb{N}}:=\left\{ u_{n}^{*}\right\} _{n\in\mathbb{N}}$
and $\left\{ y_{n}\right\} _{n\in\mathbb{N}}:=\left\{ v_{n}^{*}\right\} _{n\in\mathbb{N}}$
which are in $\ell^{q}$ and of course satisfy (\ref{eq:zHL-one})
and (\ref{eq:Equals}). This gives us that $\left\{ v_{n}^{*}\right\} _{n\in\mathbb{N}}=V_{1}\left(\left\{ u_{n}^{*}\right\} _{n\in\mathbb{N}}\right)$
for a linear operator $V_{1}$ which satisfies (\ref{eq:BothV}).
Consequently $V_{1}$ maps $E$ equipped with the equivalent norm
$\left\Vert \cdot\right\Vert _{\widetilde{E}}$ into itself with norm
not exceeding $1$ and we obtain that $\left\{ v_{n}^{*}\right\} _{n\in\mathbb{N}}\in E$
and 
\[
\left\Vert \left\{ v_{n}^{*}\right\} _{n\in\mathbb{N}}\right\Vert _{\widetilde{E}}\le\left\Vert \left\{ u_{n}^{*}\right\} _{n\in\mathbb{N}}\right\Vert _{\widetilde{E}}.
\]
In a similar way to our transition above from $\left\{ u_{n}\right\} _{n\in\mathbb{N}}$
to $\left\{ u_{n}^{*}\right\} _{n\in\mathbb{N}}$, we can now use
Lemma \ref{lem:EasyCalderon} and Remark \ref{rem:NormIsOneForEveryR}
again, this time to provide an operator $Y$ which is a norm $1$
map of $\ell^{r}$ into itself for each $r\in[1,\infty]$ and which
maps $\left\{ v_{n}^{*}\right\} _{n\in\mathbb{N}}$ to $\left\{ v_{n}\right\} _{n\in\mathbb{N}}$
and therefore shows that $\left\{ v_{n}\right\} _{n\in\mathbb{N}}\in E$
with $\left\Vert \left\{ v_{n}\right\} _{n\in\mathbb{N}}\right\Vert _{\widetilde{E}}\le\left\Vert \left\{ v_{n}^{*}\right\} _{n\in\mathbb{N}}\right\Vert _{\widetilde{E}}.$ 

Combining the preceding estimates, we see that 
\begin{align*}
\left\Vert \left\{ v_{n}\right\} _{n\in\mathbb{N}}\right\Vert _{E} & \le\left\Vert \left\{ v_{n}\right\} _{n\in\mathbb{N}}\right\Vert _{\widetilde{E}}\le\left\Vert \left\{ v_{n}^{*}\right\} _{n\in\mathbb{N}}\right\Vert _{\widetilde{E}}\le\left\Vert \left\{ u_{n}^{*}\right\} _{n\in\mathbb{N}}\right\Vert _{\widetilde{E}}\\
 & \le\left\Vert \left\{ u_{n}\right\} _{n\in\mathbb{N}}\right\Vert _{\widetilde{E}}\le C\left\Vert \left\{ u_{n}\right\} _{n\in\mathbb{N}}\right\Vert _{E}.
\end{align*}
Consequently, $E$ indeed has property $S_{q}(C)$. $\qed$

\begin{rem}
\label{rem:ArazyCwikelConjecture} Sparr obtained one version of his
main result in the larger realm of quasi-Banach spaces, i.e., for couples
of $L^{p}$ spaces where the exponents can also take positive values
less than $1$, provided certain conditions are imposed on the underlying
measure spaces. See part (ii) of Lemma 4.1 on p.~239 of \cite{SparrG1978}.
Bearing this and also Remark \ref{rem:CompareWithFSandCoauthors}
in mind one may conjecture that a version of Theorem \ref{thm:Easy}
might also hold for positive values of $p$ which are less than $1$.
As a first step towards investigating this possibility, one can attempt
to establish a variant, for quasinormed sequence spaces, of the second
part of the lemma on pp.~257-258 of \cite{ArazyJCwikelM1984}which
we used above. More explicitly, if the non-negative non-increasing
sequences $\left\{ x_{n}\right\} _{n\in\mathbb{N}}$ and $\left\{ y_{n}\right\} _{n\in\mathbb{N}}$
satisfy (\ref{eq:rtq}) does this ensure the existence of a bounded
operator, say $V_{2}:\ell^{q}\to\ell^{q}$, which, like $V_{1}$ above,
maps $\left\{ x_{n}\right\} _{n\in\mathbb{N}}$ to $\left\{ y_{n}\right\} _{n\in\mathbb{N}}$
and is also bounded from $\ell^{1}$ to itself, but can also be constructed
so that it is bounded from $\ell^{p}$ to itself for some given positive
$p$ which is less than $1$ (and perhaps even for all such $p$)?
Such a result would apparently suffice to establish a more general
version of Theorem \ref{thm:Easy} for that value of $p$ where suitably
modified definitions permit all relevant spaces to be quasinormed
instead of merely normed.
\end{rem}

\section{\label{sec:Main}Spaces with property $S_{q}(C)$ are interpolation
spaces between $\ell^{1}$ and $\ell^{q}$}

We now turn to the task of proving some sort a partial converse to
Theorem \ref{thm:Easy}, namely that sequence spaces $E$ with property
$S_{q}(C)$ should also necessarily be interpolation spaces with respect
to the couple $\left(\ell^{p},\ell^{q}\right)$ for some $p\in[1,q]$. 

It is easy to see that it is not possible to conclude that such a
sequence space $E$ is a interpolation space in general , unless some
additional hypothesis or hypotheses on $E$ are imposed. Apart from
anything else, we are obliged to require that 
\begin{equation}
E\subset\ell^{q},\label{eq:InEllQ}
\end{equation}

\noindent since any normed interpolation space with respect to $\left(\ell^{p},\ell^{q}\right)$
is necessarily contained in $\ell^{p}+\ell^{q}=\ell^{q}$. 

It seems natural to assume that $E$ is a \textit{symmetric} sequence
space, i.e., that whenever the nonincreasing rearrangements $\left\{ x_{n}^{*}\right\} _{n\in\mathbb{N}}$
and $\left\{ y_{n}^{*}\right\} _{n\in\mathbb{N}}$ of two sequences
$\left\{ x_{n}\right\} _{n\in\mathbb{N}}$ and $\left\{ y_{n}\right\} _{n\in\mathbb{N}}$
satisfy $y_{n}^{*}\le x_{n}^{*}$ and $\left\{ x_{n}\right\} _{n\in\mathbb{N}}\in E$,
then it follows that $\left\{ y{}_{n}\right\} _{n\in\mathbb{N}}\in E$
and $\left\Vert \left\{ y_{n}\right\} _{n\in\mathbb{N}}\right\Vert _{E}\le\left\Vert \left\{ x_{n}\right\} _{n\in\mathbb{N}}\right\Vert _{E}$. 

We are grateful to the authors of \cite{LeSuZa2017,LeSuZa2017B} for
a private communication in which they proved that if $E$ is any symmetric
sequence space which has property $S_{q}(C)$ then $E$ satisfies
(\ref{eq:InEllQ}). 
\begin{rem}
\label{rem:ExoticExamples} The inclusion (\ref{eq:InEllQ}) does
not follow from property $S_{q}(C)$ if $E$ is not required to be
symmetric. Given any sequence space $X$, `it is easy to construct
a sequence space $E$ which has property $S_{q}(1)$ but also contains
an isometric image of $X$ as a complemented subspace. Such constructions
also show that a normed sequence space $E$ with property $S_{q}(1)$
may fail to be contained in $\ell^{q}$ or even in $\ell^{\infty}$. 
\end{rem}
We shall show that a sequence space $E$ which has property $S_{q}(C)$
and which is contained in $\ell^{q}$ is in fact an interpolation
space with respect to $\left(\ell^{1},\ell^{q}\right)$, provided
it has one more property, a property which, at first sight, may seem
to be very close to the property that we want to prove, namely we
will be requiring that $E$ should already be an interpolation space
with respect to $\left(\ell^{1},\ell^{\infty}\right)$. But, as we
shall show in a discussion which we have deferred to Section \ref{sec:Further},
this last requirement is after all not so restrictive . It turns out
that it holds automatically (see Theorem \ref{thm:Enough} below)
if we impose another perhaps seemingly milder condition on $E$, namely
the weak version of the Fatou Property formulated above in Definition
\ref{def:WeakFatou}. But it can also hold for some spaces $E$ which
do not have the Fatou Property. (See Remark \ref{rem:NotNecessary}). 
\begin{rem}
\label{rem:NotForQuasinorms} As a digression and brief sequel to
Remark \ref{rem:ArazyCwikelConjecture} let us offer another quick
comment regarding possible attempts to extend results in this paper
to the wider context of sequence spaces which are merely quasinormed
rather than normed. For such spaces, as obvious examples show, even
the usual stronger version $WFP(1)$ of the Fatou property is not
sufficient to imply that they are (quasinormed) interpolation spaces
with respect to $\left(\ell^{1},\ell^{\infty}\right)$. (The arguments
used below in Section \ref{sec:Further}, in particular the inequality
(\ref{eq:NotForQuasinorms}) in the proof of Theorem \ref{thm:Enough},
fail for quasinormed spaces.)
\end{rem}
Here then is the main result of this section, our approximate converse
to Theorem \ref{thm:Easy}, and an essential component for the proof
of Theorem \ref{thm:ReallyMain}. 
\begin{thm}
\label{thm:Main}Let $q\in(1,\infty)$ and let $C_{1}$ and $C_{2}$
be positive constants. Let $E$ be a normed sequence space which is
a subset of $\ell^{q}$ and has property $S_{q}(C_{1})$. Suppose,
furthermore, that $E$ is a $C_{2}$-interpolation space with respect
to the couple $\left(\ell^{1},\ell^{\infty}\right)$. Then $E$ is
also a $C_{3}-K$ space and a $C_{3}$-interpolation space with respect
to the couple $\left(\ell^{1},\ell^{q}\right)$ for a constant $C_{3}$
which depends only on $q$, $C_{1}$ and $C_{2}$. 
\end{thm}
\noindent \textit{Proof.} . We will present several of the steps
of the proof as separate claims and a separate proposition. We begin
by establishing three claims about certain properties of $E$. 
\begin{claim}
\label{claim:Lattice}If $\left\{ x_{n}\right\} _{n\in\mathbb{N}}\in E$
and $\left\{ y_{n}\right\} _{n\in\mathbb{N}}$ satisfies $\left|y_{n}\right|\le\left|x_{n}\right|$
for every $n\in\mathbb{N}$, then $\left\{ y_{n}\right\} _{n\in\mathbb{N}}\in E$
and $\left\Vert \left\{ y_{n}\right\} _{n\in\mathbb{N}}\right\Vert _{E}\le C_{2}\left\Vert \left\{ x_{n}\right\} _{n\in\mathbb{N}}\right\Vert _{E}$.
\end{claim}

\begin{claim}
\label{claim:EasySqC}If $\left\{ u_{n}\right\} _{n\in\mathbb{N}}\in E$
and $\left\{ y_{n}\right\} _{n\in\mathbb{N}}$ satisfies 
\begin{equation}
\sum_{n=N}^{\infty}\left(y_{n}^{*}\right)^{q}\le\sum_{n=N}^{\infty}\left(u_{n}^{*}\right)^{q}\mbox{ for every\,}N\in\mathbb{N},\label{eq:shc}
\end{equation}
 then $\left\{ y_{n}\right\} _{n\in\mathbb{N}}\in E$ and $\left\Vert \left\{ y_{n}\right\} _{n\in\mathbb{N}}\right\Vert _{E}\le C_{1}C_{2}\left\Vert \left\{ u_{n}\right\} _{n\in\mathbb{N}}\right\Vert _{E}$.
\end{claim}

\begin{claim}
\label{claim:EisRI}If $\left\{ u_{n}\right\} _{n\in\mathbb{N}}\in E$
and $\left\{ y_{n}\right\} _{n\in\mathbb{N}}$ satisfies 
\begin{equation}
y_{n}^{*}\le u_{n}^{*}\mbox{ for every\,}n\in\mathbb{N},\label{eq:WasShc-1}
\end{equation}
 then $\left\{ y_{n}\right\} _{n\in\mathbb{N}}\in E$ and 
\begin{equation}
\left\Vert \left\{ y_{n}\right\} _{n\in\mathbb{N}}\right\Vert _{E}\le C_{2}\left\Vert \left\{ u_{n}\right\} _{n\in\mathbb{N}}\right\Vert _{E}.\label{eq:EisRI}
\end{equation}

\end{claim}
In order to prove Claim \ref{claim:Lattice} we consider the linear
map $T$ defined by 
\begin{equation}
T\left\{ h_{n}\right\} _{n\in\mathbb{N}}=\left\{ \theta_{n}h_{n}\right\} _{n\in\mathbb{N}}\label{eq:PtwsMult}
\end{equation}
where $\theta_{n}=y_{n}/x_{n}$ for all $n$ for which $x_{n}\ne0$
and $\theta_{n}=0$ whenever $x_{n}=0$. Obviously $T\in\mathcal{L}_{1}\left((\ell^{1},\ell^{\infty}\right)$
and $T\left\{ x_{n}\right\} _{n\in\mathbb{N}}=\left\{ y_{n}\right\} _{n\in\mathbb{N}}$
and so the desired conclusion follows immediately from the fact that
$E$ is a $C_{2}$-interpolation space with respect to $\left(\ell^{1},\ell^{\infty}\right)$.

In order to prove Claim \ref{claim:EasySqC}, let us first note that
the fact that $E\subset\ell^{q}$ together with (\ref{eq:shc}) ensures
that $\lim_{n\to\infty}y_{n}=0$. So there exists at least one positive
integer $n_{1}$ which has the property that $\left|y_{n}\right|\le\left|y_{n_{1}}\right|$
for every $n\in\mathbb{N}$. Having chosen a particular $n_{1}$ with
this property, let $\left\{ z_{n}\right\} _{n\in\mathbb{N}}$ be a
sequence which satisfies 
\begin{equation}
z_{n_{1}}\ge\left|y_{n_{1}}\right|\label{eq:zyno}
\end{equation}
and $z_{n}=\left|y_{n}\right|$ for every $n\in\mathbb{N}\setminus\left\{ n_{1}\right\} $.
Then the nonincreasing rearrangements $\left\{ z_{n}^{*}\right\} _{n\in\mathbb{N}}$
and $\left\{ y_{n}^{*}\right\} _{n\in\mathbb{N}}$ of $\left\{ z_{n}\right\} _{n\in\mathbb{N}}$
and $\left\{ y_{n}\right\} _{n\in\mathbb{N}}$ satisfy $z_{1}^{*}=z_{n_{1}}\ge\left|y_{n_{1}}\right|=y_{1}^{*}$
and $z_{n}^{*}=y_{n}^{*}$ for every $n\ge2$. In view of (\ref{eq:shc}),\textit{
every} sequence $\left\{ z_{n}\right\} _{n\in\mathbb{N}}$ which is
obtained from $\left\{ y_{n}\right\} _{n\in\mathbb{N}}$ in the particular
way specified just above satisfies 
\begin{equation}
\sum_{n=N}^{\infty}\left(z_{n}^{*}\right)^{q}\le\sum_{n=N}^{\infty}\left(u_{n}^{*}\right)^{q}\label{eq:znx}
\end{equation}
for every $N\ge2$. Since (\ref{eq:shc}) holds also for $N=1$ it
is clear that we can choose a particular value for $z_{n_{1}}$ which
will satisfy (\ref{eq:zyno}) and will also satisfy 
\begin{equation}
\sum_{n=1}^{\infty}\left(z_{n}^{*}\right)^{q}=\sum_{n=1}^{\infty}\left(u_{n}^{*}\right)^{q}.\label{eq:jomp}
\end{equation}
We know that the sum of the series in (\ref{eq:jomp}) is finite,
since we have $\left\{ u_{n}\right\} _{n\in\mathbb{N}}\in\ell^{q}$.
So we can now use the simple ``reverse'' of the simple argument
that appeared at the beginning of the proof of Theorem \ref{thm:Easy}.
We obtain, from (\ref{eq:jomp}) and (\ref{eq:znx}), that, for each
$N\in\mathbb{N}$, 
\begin{align*}
\sum_{m=1}^{N}\left(u_{m}^{*}\right)^{q} & =\left(\sum_{m=1}^{\infty}\left(u_{m}^{*}\right)^{q}-\sum_{m=N+1}^{\infty}\left(u_{m}^{*}\right)^{q}\right)\\
 & =\sum_{m=1}^{\infty}\left(z_{m}^{*}\right)^{q}-\sum_{m=N+1}^{\infty}\left(u_{m}^{*}\right)^{q}\\
 & \le\sum_{m=1}^{\infty}\left(z_{m}^{*}\right)^{q}-\sum_{m=N+1}^{\infty}\left(z_{m}^{*}\right)^{q}=\sum_{m=1}^{N}\left(z_{m}^{*}\right)^{q}
\end{align*}
 This last inequality for each $N$, together with (\ref{eq:znx})
and the facts that the sequence $\left\{ u_{n}\right\} _{n\in\mathbb{N}}$
is an element of $E$ and $E$ has property $S_{q}(C_{1})$ all combine
to imply that $\left\{ z_{n}\right\} _{n\in\mathbb{N}}$ is an element
of $E$ and satisfies \textbf{
\[
\left\Vert \left\{ z_{n}\right\} _{n\in\mathbb{N}}\right\Vert _{E}\le C_{1}\left\Vert \left\{ u_{n}\right\} _{n\in\mathbb{N}}\right\Vert _{E}.
\]
}So we can now also use Claim \ref{claim:Lattice} to obtain that
$\left\{ y_{n}\right\} _{n\in\mathbb{N}}$ is also an element of $E$
and 
\[
\left\Vert \left\{ y_{n}\right\} _{n\in\mathbb{N}}\right\Vert _{E}\le C_{2}\left\Vert \left\{ z_{n}\right\} _{n\in\mathbb{N}}\right\Vert _{E}\le C_{1}C_{2}\left\Vert \left\{ u_{n}\right\} _{n\in\mathbb{N}}\right\Vert _{E}
\]
and this of course completes the proof of Claim \ref{claim:EasySqC}.

With regard to Claim \ref{claim:EisRI}, let us first observe in passing
that Claim \ref{claim:EasySqC} is immediately applicable and gives
us that $\left\{ y_{n}\right\} _{n\in\mathbb{N}}\in E$. It also gives
us the inequality $\left\Vert \left\{ y_{n}\right\} _{n\in\mathbb{N}}\right\Vert _{E}\le C_{1}C_{2}\left\Vert \left\{ u_{n}\right\} _{n\in\mathbb{N}}\right\Vert _{E}$
which may be weaker and certainly cannot (cf.~Remark \ref{rem:BiggerThanOne})
be stronger than (\ref{eq:EisRI}). In order to obtain (\ref{eq:EisRI})
we use Lemma \ref{lem:EasyCalderon} to supply us with $W$ and $Y$
in $\mathcal{L}_{1}\left(\left(\ell^{1},\ell^{\infty}\right)\right)$
for which $W\left(\left\{ u_{n}\right\} _{n\in\mathbb{N}}\right)=\left\{ u_{n}^{*}\right\} _{n\in\mathbb{N}}$
and $Y\left(\left\{ y_{n}^{*}\right\} _{n\in\mathbb{N}}\right)=\left\{ y_{n}\right\} _{n\in\mathbb{N}}$.
Then, in view of (\ref{eq:WasShc-1}), there is another obvious map
$T\in\mathcal{L}_{1}\left(\left(\ell^{1},\ell^{\infty}\right)\right)$
of pointwise multiplication (cf.~(\ref{eq:PtwsMult})) which satisfies
$T\left(\left\{ u_{n}^{*}\right\} _{n\in\mathbb{N}}\right)=\left\{ y_{n}^{*}\right\} _{n\in\mathbb{N}}$.
The composed operator $YTW$ is also of course in $\mathcal{L}_{1}\left(\left(\ell^{1},\ell^{\infty}\right)\right)$
and maps $\left\{ u_{n}\right\} _{n\in\mathbb{N}}$ to $\left\{ y_{n}\right\} _{n\in\mathbb{N}}$
immediately giving us (\ref{eq:EisRI}) and completing the proof of
Claim \ref{claim:EisRI}.

\medskip{}

Now let us describe the main task which needs to be performed in order
to prove Theorem \ref{thm:Main}. In view of Fact \ref{fac:KspaceIsInterp},
it will suffice to show that $E$ is a $C_{3}-K$ space for the couple
$\left(\ell^{1},\ell^{q}\right)$ and for some constant $C_{3}$ which
depends only on $q$, $C_{1}$ and $C_{2}$. Accordingly, let us fix
an arbitrary sequence $x=\left\{ x_{n}\right\} _{n\in\mathbb{N}}$
in $E$ and let $y=\left\{ y_{n}\right\} _{n\in\mathbb{N}}$ be an
arbitrary sequence in $\ell^{q}$ which satisfies 
\begin{equation}
K(t,y;\ell^{1},\ell^{q})\le K(t,x;\ell^{1},\ell^{q})\mbox{\,\ for every\,}t>0,\label{eq:XandY}
\end{equation}
and we will, after a number of intermediate steps, complete the proof
of this theorem by showing that 
\begin{equation}
y\in E\mbox{\,\ and\,}\left\Vert y\right\Vert _{E}\le C_{3}\left\Vert x\right\Vert _{E}\label{eq:toep}
\end{equation}
where $C_{3}$ is as specified above.

We shall use a result about the $K$-functional for $\left(L^{1},L^{q}\right)$
(and therefore ultimately also for $\left(\ell^{1},\ell^{q}\right)$)
which is a special case of a result obtained by Tord Holmstedt. (See
Theorem 4.1 on p.~189 of \cite{HolmstedtT1979}.) Holmstedt's result
ensures that, for some suitable constant $C(q)$, depending only on
$q$, for $\alpha=q/(q-1)$ and for every $f\in L^{1}+L^{q}$ and
for every $t>0$, we have 
\begin{equation}
K(t,f;L^{1},L^{q})\le\int_{0}^{t^{\alpha}}f^{*}(s)ds+t\left(\int_{t^{\alpha}}^{\infty}\left(f^{*}(s)\right)^{q}ds\right)^{1/q}\le C(q)K(t,f;L^{1},L^{q}).\label{eq:Holmstedt}
\end{equation}

\begin{rem}
If we substitute $p=1$ in an inequality mentioned on p.~255 of \cite{ArazyJCwikelM1984}
we can apparently obtain that 
\[
C(q)\le\max\left\{ \left(1+2(q-1)^{-1/q}\right)\,,\,\left(2^{1/q}+(1-1/q)^{-1/q}\right)\right\} .
\]
We also note that an exact formula for $K(t,f;L^{1},L^{q})$ and a
precise description of the two functions whose sum is $f$ and for
which the infimum for calculating $K(t,f;L^{1},L^{q})$ is attained
can be found in \cite[Theorem 1', p.~324 and Remark 3, p.~325]{NilssonPPeetreJ1986}.
But at this stage we do not see any way of using these results to
obtain a version of our theorem with a better estimate for $C_{3}$.
\end{rem}
In fact Holmstedt establishes his result for couples of $L^{p}$ spaces
on an arbitrary measure space, thus also including the case where
the measure space is $\mathbb{N}$ equipped with counting measure.
Because of this, or via an easy argument using (\ref{eq:Holmstedt})
only in the case where the measure space is $\left(0,\infty\right)$
with Lebesgue measure, we can readily obtain that, for any sequence
$u=\left\{ u_{n}\right\} _{n\in\mathbb{N}}$ in $\ell^{p}$ whose
nonincreasing rearrangement is the sequence $\left\{ u_{n}^{*}\right\} _{n\in\mathbb{N}}$,
we have
\begin{equation}
K(t,u;\ell^{1},\ell^{q})\le\int_{0}^{t^{\alpha}}\varphi(s)ds+t\left(\int_{t^{\alpha}}^{\infty}\left(\varphi(s)\right)^{q}ds\right)^{1/q}\le C(q)K(t,u;\ell^{1},\ell^{q})\label{eq:SeqHolm}
\end{equation}
for all $t>0$, where $\varphi:[0,\infty)\to[0,\infty)$ is the nonincreasing
function $\varphi=\sum_{n=1}^{\infty}u_{n}^{*}\chi_{[n-1,n)}$.

Let the sequences $\left\{ x_{n}^{*}\right\} _{n\in\mathbb{N}}$ and
$\left\{ y_{n}^{*}\right\} _{n\in\mathbb{N}}$ be the nonincreasing
rearrangements of the special (but arbitrary) sequences $x=\left\{ x_{n}\right\} _{n\in\mathbb{N}}$
and $y=\left\{ y_{n}\right\} _{n\in\mathbb{N}}$ which were chosen
above. Since $x$ and $y$ satisfy (\ref{eq:XandY}), it follows,
also from (\ref{eq:SeqHolm}), that the nonincreasing functions $g$
and $h$ defined by 
\begin{equation}
g=\sum_{n=1}^{\infty}y_{n}^{*}\chi_{[n-1,n)}\mbox{\,\,\ and\,\,}h=C(q)\sum_{n=1}^{\infty}x_{n}^{*}\chi_{[n-1,n)}\label{eq:defgh}
\end{equation}
 satisfy 
\begin{align}
 & \int_{0}^{t^{\alpha}}g(s)ds+t\left(\int_{t^{\alpha}}^{\infty}\left(g(s)\right)^{q}ds\right)^{1/q}\label{eq:firstLine}\\
\le & \, C(q)K(t,y;\ell^{1},\ell^{q})\nonumber \\
\le & \, C(q)K(t,x;\ell^{1},\ell^{q})=K(t,C(q)x;\ell^{1},\ell^{q})\nonumber \\
\le & \int_{0}^{t^{\alpha}}h(s)ds+t\left(\int_{t^{\alpha}}^{\infty}\left(h(s)\right)^{q}ds\right)^{1/q}\nonumber 
\end{align}
for all $t>0$. Let $\varepsilon$ be an arbitrary positive number
and set 
\begin{equation}
f=(1+\varepsilon)h.\label{eq:deff}
\end{equation}
We can of course assume that $y$ is not identically zero and that,
consequently, the expression in the first line (\ref{eq:firstLine})
of the preceding calculation is strictly positive for all $t>0$.
Therefore that preceding calculation implies that 
\[
\int_{0}^{t^{\alpha}}g(s)ds+t\left(\int_{t^{\alpha}}^{\infty}\left(g(s)\right)^{q}ds\right)^{1/q}<\int_{0}^{t^{\alpha}}f(s)ds+t\left(\int_{t^{\alpha}}^{\infty}\left(f(s)\right)^{q}ds\right)^{1/q}
\]
for every $t>0$. This in turn means that, at each point $t\in(0,\infty)$,
at least one of the two inequalities 
\begin{equation}
\int_{0}^{t}g(s)ds<\int_{0}^{t}f(s)ds\label{eq:FirstIneq}
\end{equation}
and 
\begin{equation}
\int_{t}^{\infty}\left(g(s)\right)^{q}ds<\int_{t}^{\infty}\left(f(s)\right)^{q}ds\label{eq:SecondIneq}
\end{equation}
must hold. Since $\left\Vert y\right\Vert _{\ell^{q}}=\left\Vert g\right\Vert _{L^{q}}$
and we are assuming that $y\ne0$, we see that (\ref{eq:SecondIneq})
must hold also at $t=0$.

Let $A$ be the set of all points in $[0,\infty)$ which satisfy (\ref{eq:FirstIneq}).
Let $B$ be the set of all points in $[0,\infty)$ which satisfy (\ref{eq:SecondIneq}).
Our remarks just above when we introduced the inequalities (\ref{eq:FirstIneq})
and (\ref{eq:SecondIneq}) can be equivalently restated as the formula
\begin{equation}
A\cup B=[0,\infty).\label{eq:AuB}
\end{equation}
together with the fact that $0\in B$. 

(Any reader who happens to be familiar with one of the proofs in \cite{CwikelM1976}
and later in \cite{ArazyJCwikelM1984} will recognize some considerable
similarities (and also observe some considerable differences) between
the general approach in those proofs and the path that we have begun
following here and will continue following now.) 

Each of the four integrals which appear in (\ref{eq:FirstIneq}) and
(\ref{eq:SecondIneq}) is a continuous function of $t$. Therefore
the intersections of $A$ and of $B$ with $\left(0,\infty\right)$
must both be open sets.

It will be evident from later considerations that the proof of our
theorem would be considerably simpler if we could assert that each
of the sets $A$ and $B\setminus\{0\}$ is the union of a finite or
infinite sequence of open intervals all of whose finite endpoints
are integers. We cannot assert this, but the following two claims
will make it possible for us to, in some sense, ``almost'' assert
this. (Cf.~Proposition \ref{prop:FourSequences} below.) The first
of these claims deals with the set $A$.
\begin{claim}
\label{claim:Astuff}Suppose that $0\le a<\widetilde{a}\le\infty$
and that the open interval $\left(a,\widetilde{a}\right)$ is contained
in $A$ and that $a$ is a boundary point of $A$. Suppose also that,
if $\widetilde{a}<\infty$, then it is also a boundary ponit of $A$.
Let $a^{\clubsuit}$ be the largest integer which satisfies $a^{\clubsuit}\le a$.
Then 
\begin{equation}
\int_{a^{\clubsuit}}^{t}g(s)ds\le\int_{a^{\clubsuit}}^{t}f(s)ds\mbox{ for all\,}t\in\left[a^{\clubsuit},\max\left\{ \widetilde{a},a^{\clubsuit}+1\right\} \right]\label{eq:gfaa}
\end{equation}
and 
\begin{equation}
g(a^{\clubsuit})\le f(a^{\clubsuit}).\label{eq:gfaaN}
\end{equation}
If $\widetilde{a}<\infty$ and $a^{\diamondsuit}$ is the largest
integer which satisfies $a^{\diamondsuit}\le\max\left\{ \widetilde{a},a^{\clubsuit}+1\right\} $
then at least one of the two integers $a^{\diamondsuit}$ and $a^{\diamondsuit}+1$
must be an element of $B\setminus A$.
\end{claim}
\textit{Proof of Claim \ref{claim:Astuff}.} Let us first prove a
possibly weaker version of (\ref{eq:gfaa}) namely, that 
\begin{equation}
\int_{a^{\clubsuit}}^{t}g(s)ds\le\int_{a^{\clubsuit}}^{t}f(s)ds\mbox{ for all\,}t\in\left[a^{\clubsuit},\widetilde{a}\right].\label{eq:gfaaWeak}
\end{equation}
The fact that $a\in\partial A$ implies that 
\[
\int_{0}^{a}g(s)ds=\int_{0}^{a}f(s)ds.
\]
Therefore, for each $t\in(a,\widetilde{a})$ we have 
\begin{align}
 & {\phantom{=}}\int_{a}^{t}f(s)ds-\int_{a}^{t}g(s)ds\nonumber \\
 & =\int_{0}^{t}f(s)ds-\int_{0}^{t}g(s)ds>0.\label{eq:ut}
\end{align}
We shall consider two cases, where$a$ is, or is not, an integer.

If $a$ is an integer then $a=a^{\clubsuit}$ and so (\ref{eq:ut})
shows that $\int_{a^{\clubsuit}}^{t}g(s)ds\le\int_{a^{\clubsuit}}^{t}f(s)ds$
for all $t\in(a,\widetilde{a})=(a^{\clubsuit},\widetilde{a})$ and
therefore also, by continuity, for $t=\widetilde{a}$ (even if $\widetilde{a}=\infty)$
and even if $\int_{a^{\clubsuit}}^{\infty}g(s)ds$ or $\int_{a^{\clubsuit}}^{\infty}f(s)ds$
is infinite). It also holds (trivially!) for $t=a^{\clubsuit}$ so
this completes the proof of (\ref{eq:gfaaWeak}) in this case. 

In the remaining case, where $a$ is not an integer, the intervals
$\left[a^{\clubsuit},a\right]$ and $[a,a^{\clubsuit}+1)$ both have
positive length and the function $f-g$ takes the same constant value,
$r:=f(a^{\clubsuit})-g(a^{\clubsuit})$, on both of them. Since $\int_{a}^{t}f(s)ds-\int_{a}^{t}g(s)ds>0$
for every $t\in\left(a,\min\left\{ \widetilde{a},a^{\clubsuit}+1\right\} \right)$
we see that $r>0$. This will now enable us to prove (\ref{eq:gfaaWeak})
also in this case:

First, for each $t\in[a^{\clubsuit},a]$, we have 
\[
\int_{a^{\clubsuit}}^{t}f(s)ds-\int_{a^{\clubsuit}}^{t}g(s)ds=(t-a^{\clubsuit})r\ge0.
\]
Then, for each $t\in[a,\widetilde{a})$ we have, also using (\ref{eq:ut}),
that
\begin{align*}
 & {\phantom{=}}\int_{a^{\clubsuit}}^{t}f(s)ds-\int_{a^{\clubsuit}}^{t}g(s)ds\\
 & =\int_{a}^{t}f(s)ds-\int_{a}^{t}g(s)ds+(a-a^{\clubsuit})r\\
 & >\int_{a}^{t}f(s)ds-\int_{a}^{t}g(s)ds\ge0.
\end{align*}
So $\int_{a^{\clubsuit}}^{t}g(s)ds\le\int_{a^{\clubsuit}}^{t}g(s)ds$
holds for all $t\in[a^{\clubsuit},\widetilde{a})$ and we can also
deduce that it holds for $t=\widetilde{a}$ by passing to the limit
as $t\to\widetilde{a}$, whether or not $\widetilde{a}$ is finite.
This completes the proof of (\ref{eq:gfaaWeak}) in the second case
and therefore in general.

Before proving (\ref{eq:gfaa}) we shall prove (\ref{eq:gfaaN}).
In the case where $a$ is not an integer, our previous reasoning has
already shown that $f(a^{\clubsuit})-g(a^{\clubsuit})>0$ and thus
established (\ref{eq:gfaaN}) (even with strict inequality). Regardless
of that, we can now, whether or not $a$ is an integer, obtain (\ref{eq:gfaaN})
from (\ref{eq:gfaaWeak}) by choosing $t$ (as before) such that $a^{\clubsuit}<t<\min\left\{ \widetilde{a},a^{\clubsuit}+1\right\} $
and dividing by $t-a^{\clubsuit}$. This gives us that $\frac{1}{t-a^{\clubsuit}}\int_{a^{\clubsuit}}^{t}g(s)ds\le\frac{1}{t-a^{\clubsuit}}\int_{a^{\clubsuit}}^{t}f(s)ds$
which, is exactly (\ref{eq:gfaaN}). 

It is an obvious consequence of (\ref{eq:gfaaN}) that $\int_{a^{\clubsuit}}^{t}g(s)ds\le\int_{a^{\clubsuit}}^{t}g(s)ds$
also holds for all $t$ in the interval $\left[a^{\clubsuit},a^{\clubsuit}+1\right]$
and so, by (\ref{eq:gfaaWeak}), it must also hold for all $t\in\left[a^{\clubsuit},\max\left\{ \widetilde{a},a^{\clubsuit}+1\right\} \right]$
and this completes the proof of (\ref{eq:gfaa}). 

It remains to prove the claim concerning $a^{\diamondsuit}$. Here
we are supposing that $\widetilde{a}$ is finite. Therefore, as stipulated
in the statement of Claim \ref{claim:Astuff}, we must have $\widetilde{a}\in\partial A$.
This implies that the continuous function $\varphi(t):=\int_{0}^{t}f(s)-g(s)ds$
must vanish at $t=\widetilde{a}$. Since $\left(a,\widetilde{a}\right)\subset A$
we must also have $\varphi(t)>0$ for every $t\in(a,\widetilde{a})$. 

Suppose that $\widetilde{a}\in(n-1,n)$ for some $n\in\mathbb{N}$.
Since $\varphi$ is affine on $\left[n-1,n\right]$ and vanishes at
$\widetilde{a}$ and is strictly positive at some, and therefore all
points in $[n-1,\widetilde{a})$ we must also have $\varphi(t)<0$
and therefore $t\notin A$ for all $t\in(\widetilde{a},n]$. In particular
we must have that $n\notin A$ and therefore very obviously (recalling
(\ref{eq:AuB})) that in fact $n\in B\setminus A$.

We need to apply the preceding reasoning in two separate cases:

First, if $\widetilde{a}<a^{\clubsuit}+1$, we must have $a^{\diamondsuit}=a^{\clubsuit}+1$
and so 
\[
a^{\diamondsuit}-1=a^{\clubsuit}\le a<\widetilde{a}<a^{\clubsuit}+1=a^{\diamondsuit}.
\]
So $\widetilde{a}\in(a^{\diamondsuit}-1,a^{\diamondsuit})$ and the
above reasoning gives us that $a^{\diamondsuit}\in B\setminus A$.
Alternatively, if $\widetilde{a}\ge a^{\clubsuit}+1$, we must have
$a^{\diamondsuit}\le\widetilde{a}<a^{\diamondsuit}+1$. If $\widetilde{a}=a^{\diamondsuit}$
then, since $\widetilde{a}\in\partial A\subset[0,\infty)\setminus A$
we must again have $a^{\diamondsuit}\in B\setminus A$. Otherwise
$\widetilde{a}$ is an interior point of $\left[a^{\diamondsuit},a^{\diamondsuit}+1\right]$
and the above reasoning using the function $\varphi$ gives us in
this case that $a^{\diamondsuit}+1\in B\setminus A$. 

\smallskip{}

We now turn to a claim about the set $B$ which is quite closely analogous
to a weaker version of Claim \ref{claim:Astuff}. 
\begin{rem}
\label{rem:AboutStrongerVersion}Apparently a stronger version of
this claim with conclusions analogous to latter parts of Claim \ref{claim:Astuff}
also holds. We will not formulate it here, but will roughly indicate
one part of such a stronger version below in our proof of (\ref{eq:fgbheartz}).\end{rem}
\begin{claim}
\label{claim:Bstuff}Suppose that $0\le\beta<b<\infty$ and that the
open interval $\left(\beta,b\right)$ is contained in $B$ and that
$b$ is a boundary point of $B$. Let $b^{\diamondsuit}$ be the smallest
integer which satisfies $b\le b^{\diamondsuit}$. Then
\begin{equation}
\int_{t}^{b^{\diamondsuit}}\left(g(s)\right)^{q}ds\le\int_{t}^{b^{\diamondsuit}}\left(f(s)\right)^{q}ds\mbox{ for all\,}t\in[\beta,b^{\diamondsuit}].\label{eq:ewb}
\end{equation}

\end{claim}
\textit{Proof of Claim \ref{claim:Bstuff}.} As the reader will probably
notice, this proof is an almost obvious ``transformation'' of the
first part of the proof of Claim \ref{claim:Astuff}, the part from
which we obtained (\ref{eq:gfaaWeak}). Nevertheless we have chosen
to present it explicitly. 

The fact that $b\in\partial B$ implies that
\[
\int_{b}^{\infty}\left(g(s)\right)^{q}ds=\int_{b}^{\infty}\left(f(s)\right)^{q}ds.
\]
Therefore, for each $t\in(\beta,b)$, we have 
\begin{align}
 & {\phantom{=}}\int_{t}^{b}\left(f(s)\right)^{q}ds-\int_{t}^{b}\left(g(s)\right)^{q}ds\nonumber \\
 & =\int_{t}^{\infty}\left(f(s)\right)^{q}ds-\int_{t}^{\infty}\left(g(s)\right)^{q}ds>0.\label{eq:put}
\end{align}
If $b\in\mathbb{N}$ so that $b=b^{\diamondsuit}$ this gives us (\ref{eq:ewb})
for all $t\in(\beta,b)=(\beta,b^{\diamondsuit})$ and therefore also,
by continuity, for $t=\beta$ and (obviously!) for $t=b^{\diamondsuit}$
and completes the proof of this claim in this case. 

Otherwise, if $b$ is not an integer, then $f^{p}-g^{p}$ takes the
same constant value, say $r$, on both of the intervals $[b^{\diamondsuit}-1,b]$
and $[b,b^{\diamondsuit})$. Since $\int_{t}^{b}\left(f(s)\right)^{q}ds-\int_{t}^{b}\left(g(s)\right)^{q}ds>0$
for every $t\in(\max\left\{ \beta,b^{\diamondsuit}-1\right\} ,b)$
we see that $r>0$. This will now enable us to establish (\ref{eq:ewb})
for all $t\in[\beta,b^{\diamondsuit}]$.

First, for each $t\in[b,b^{\diamondsuit}]$, we have 
\[
\int_{t}^{b^{\diamondsuit}}\left(f(s)\right)^{q}ds-\int_{t}^{b^{\diamondsuit}}\left(g(s)\right)^{q}ds=(b^{\diamondsuit}-t)r\ge0.
\]
Then, for each $t\in[\beta,b]$ we have, also using (\ref{eq:put}),
that 

\begin{align*}
 & {\phantom{=}}\int_{t}^{b^{\diamondsuit}}\left(f(s)\right)^{q}ds-\int_{t}^{b^{\diamondsuit}}\left(g(s)\right)^{q}ds\\
 & =\int_{t}^{b}\left(f(s)\right)^{q}ds-\int_{t}^{b}\left(g(s)\right)^{q}ds+(b^{\diamondsuit}-b)r\\
 & >\int_{t}^{b}\left(f(s)\right)^{q}ds-\int_{t}^{b}\left(g(s)\right)^{q}ds\ge0.
\end{align*}
So (\ref{eq:ewb}) indeed holds for all $t\in[\beta,b^{\diamondsuit}]$
also when $b$ is not an integer. This completes the proof of Claim
\ref{claim:Bstuff}. 

\medskip{}

Now that we are equipped with the preceding two claims, we can return
to our main task of showing that the two particular sequences $x$
and $y$ which we introduced just before (\ref{eq:XandY}) indeed
satisfy (\ref{eq:toep}). 

Our main step for doing this will be expressed by the following proposition:

\begin{prop}
\label{prop:FourSequences}There exist four special sequences $\left\{ A_{n}\right\} _{n\in\mathbb{N}}$,
$\left\{ B_{n}\right\} _{n\in\mathbb{N}}$, $\left\{ \Omega_{n}\right\} _{n\in\mathbb{N}}$
and $\left\{ \Gamma_{n}\right\} _{n\in\mathbb{N}}$ of subsets of
$[0,\infty)$ which satisfy 
\begin{equation}
[0,\infty)=\bigcup_{n\in\mathbb{N}}\left(A_{n}\cup B_{n}\cup\Gamma_{n}\right)\label{eq:UnionIsAll}
\end{equation}

and also have the following special properties for every $n\in\mathbb{N}$:

(i) Each of the sets $A_{n}$, $B_{n}$, $\Omega_{n}$ and $\Gamma_{n}$
is either the empty set or is a semi-open interval of the form $[\gamma,\delta)$
where $\gamma$ is a non-negative integer and $\delta$ is either
an integer or equals $\infty$.

(ii) $A_{n}<A_{n+1}$, $B_{n}<B_{n+1}$, $\Omega_{n}<\Omega_{n+1}$
and $\Gamma_{n}<\Gamma_{n+1}$, where, in general, the notation $G<H$
means that if $s\in G$ and $t\in H$ then $s<t$. (So of course the
empty set $\emptyset$ satisfies $\emptyset<G$ and $G<\emptyset$
for any subset $G$ of $[0,\infty)$.)

(iii) $\Omega_{n}<\Gamma_{n}$.

(iv) The special functions $f$ and $g$ introduced above satisfy
\begin{equation}
\int_{0}^{t}g(s)\chi_{A_{n}}(s)ds\le\int_{0}^{t}f(s)\chi_{A_{n}}(s)ds\mbox{ for all\,\ \ensuremath{t\ge0}}\label{eq:A-ineq}
\end{equation}
and 
\begin{equation}
\int_{t}^{\infty}\left(g(s)\right)^{q}\chi_{B_{n}}(s)ds\le\int_{t}^{\infty}\left(f(s)\right)^{q}\chi_{B_{n}}(s)ds\mbox{ for all\,\ \ensuremath{t\ge0.}}\label{eq:B-ineq}
\end{equation}
(These inequalities are of course trivially true if $A_{n}$ or $B_{n}$
is empty.)

(v) The set $\Omega_{n}$ is non-empty if and only if $\Gamma_{n}$
is non-empty. When these sets are non-empty they both have length
$1$. The constant value of $g$ on $\Gamma_{n}$ does not exceed
the constant value of $f$ on $\Omega_{n}$. 
\end{prop}
\smallskip{}

Before we describe how to construct these four sequences we should
mention that the intervals $A_{n}$ are not necessarily contained
in $A$ and the intervals $B_{n}$ are not necessarily contained in
$B$. However $A_{n}$ and $B_{n}$, when they are non-empty, have
special connections with $A$ and with $B$ respectively.

We are going to prove Proposition \ref{prop:FourSequences} with the
help of a special procedure. Let us refer to it as ``\textit{Procedure
P}''.

The starting point of Procedure P will always be an integer $n\in\mathbb{N}$
which can be considered to be simply a ``label'', and another integer
in $B\setminus A$ which we will denote by $b_{n}^{\clubsuit}$. We
will use $b_{n}^{\clubsuit}$ to construct four sets $A_{n}$, $B_{n}$,
$\Omega_{n}$ and $\Gamma_{n}$ with several special properties, among
them the property that 
\begin{equation}
\mbox{\ensuremath{A_{n}}, \ensuremath{B_{n},}\ensuremath{\Omega_{n}}and \ensuremath{\Gamma_{n}}are all contained in }[b_{n}^{\clubsuit},\infty).\label{eq:AllContained}
\end{equation}
The set $B_{n}$ will be the interval $[b_{n}^{\clubsuit},b_{n}^{\diamondsuit})$
where $b_{n}^{\diamondsuit}$ is either $\infty$ or a certain integer
strictly greater than $b_{n}^{\clubsuit}$. The set $A_{n}$ will
either be empty or will be the interval $[a_{n}^{\clubsuit},a_{n}^{\diamondsuit})$
where $a_{n}^{\clubsuit}$ is a certain integer satisfying $b_{n}^{\clubsuit}\le a_{n}^{\clubsuit}$
and $a_{n}^{\diamondsuit}$ is either $\infty$ or is another integer
strictly larger than $a_{n}^{\clubsuit}$. When $A_{n}$ is non empty
and $a_{n}^{\diamondsuit}<\infty$ we will have $\Omega_{n}=[a_{n}^{\clubsuit},a_{n}^{\clubsuit}+1)$
and $\Gamma_{n}=[a_{n}^{\diamondsuit},a_{n}^{\diamondsuit}+1)$. Otherwise
$\Omega_{n}$ and $\Gamma_{n}$ will both be empty. 

We will show that these four sets satisfy, for this particular choice
of $n$, all the properties listed in paragraphs (i), (iii), (iv)
and (v) of Proposition \ref{prop:FourSequences}. We will deal with
the properties in paragraph (ii) at a later stage.

As already hinted above, the outcome of Procedure P can be different
for different choices of the integer $b_{n}^{\clubsuit}$. We will
see that in fact there are three possible (mutually exclusive) outcomes.
We shall label them as (O-1), (O-2) and (O-3) and describe them now:

\textbf{(O-1)} The sets $A_{n}$, $B_{n}$, $\Omega_{n}$ and $\Gamma_{n}$
are all bounded and non-empty and satisfy
\begin{equation}
[b_{n}^{\clubsuit},b_{n+1}^{\clubsuit})\subset A_{n}\cup B_{n}\cup\Omega_{n}\cup\Gamma_{n}\subset[b{}_{n}^{\clubsuit},b_{n+1}^{\clubsuit}+1)\label{eq:ABGclub}
\end{equation}
for a certain integer $b_{n+1}^{\clubsuit}$ in $B\setminus A$ which
satisfies 
\begin{equation}
b_{n+1}^{\clubsuit}\ge b_{n}^{\clubsuit}+1.\label{eq:bnClubGrows}
\end{equation}

\textbf{(O-2)} The set $B_{n}$ coincides with $[b_{n}^{\clubsuit},\infty)$
and $A_{n}=\Omega_{n}=\Gamma_{n}=\emptyset$.

\textbf{(O-3)} The set $A_{n}\cup B_{n}$ coincides with $[b_{n}^{\clubsuit},\infty)$
and $\Omega_{n}=\Gamma_{n}=\emptyset$.

The reader can probably already guess how we are going to use appropriate
iterations of Procedure P to complete the proof of Proposition \ref{prop:FourSequences}.
But we will explain that precisely afterwards.

Now we are ready to begin explicitly describing Procedure P. As already
mentioned above, it starts after we choose some $n\in\mathbb{N}$
and a non-negative integer $b_{n}^{\clubsuit}\in B\setminus A$. (These
might have been generated by a previous application of Procedure P.) 

Let us first describe what happens if our integer $b_{n}^{\clubsuit}$
has the property that 
\begin{equation}
[b_{n}^{\clubsuit},\infty)\subset B.\label{eq:GivesO2}
\end{equation}
In this case we simply choose $B_{n}=[b_{n}^{\clubsuit},\infty)$
and $A_{n}=\Omega_{n}=\Gamma_{n}=\emptyset$ which completes the procedure.
I.e., in this case we have obtained outcome \textbf{(O-2)} from among
the three possible outcomes mentioned above and we have chosen $b_{n}^{\diamondsuit}$
to be $\infty$. In this case, in view of the definition of $B$,
we obtain that (\ref{eq:B-ineq}) holds for this particular set $B_{n}$
for this particular choice of $n$. It is also obvious that all other
properties listed in paragraphs (i), (iii), (iv) and (v) hold for
this particular value of $n$ and the sets $A_{n}$, $B_{n}$, $\Omega_{n}$
and $\Gamma_{n}$ that we have associated with $n$. Most of these
properties are trivialities because the relevant sets are empty. We
also obviously have (\ref{eq:AllContained}).

Let us next describe Procedure P and its consequences in the remaining
case where (\ref{eq:GivesO2}) does not hold, i.e., when $[b_{n}^{\clubsuit},\infty)$
contains at least one point which is not in $B$. In this case, since
$t\mapsto\int_{t}^{\infty}\left(f(s)\right)^{q}-\left(g(s)\right)^{q}ds$
is a continuous function which is strictly positive at $t=b_{n}^{\clubsuit}$,
we see that the supremum $b_{n}:=\sup\left\{ s\ge b_{n}^{\clubsuit}:[b_{n}^{\clubsuit},s]\subset B\right\} $
must be strictly positive but finite, and that $[b_{n}^{\clubsuit},b_{n})\subset B$
and $b_{n}\in\partial B$. All this enables us to apply Claim \ref{claim:Bstuff},
for the particular choices $\beta=b_{n}^{\clubsuit}$ and $b=b_{n}$
and thereby to obtain that 
\begin{equation}
\int_{t}^{b_{n}^{\diamondsuit}}\left(g(s)\right)^{q}ds\le\int_{t}^{b_{n}^{\diamondsuit}}\left(f(s)\right)^{q}ds\mbox{ for all\,}t\in[b_{n}^{\clubsuit},b_{n}^{\diamondsuit}]\label{eq:getheart}
\end{equation}
where, in this case we define $b_{n}^{\diamondsuit}$ the smallest
integer which satisfies $b_{n}\le b_{n}^{\diamondsuit}$. We then
set $B_{n}=[b_{n}^{\clubsuit},b_{n}^{\diamondsuit})$. This ensures
that (\ref{eq:B-ineq}) holds for this choice of $n$ and this choice
of $B_{n}$ since, obviously, it is exactly the same as (\ref{eq:getheart}). 

Since $b_{n}$ is not in $B$, it must be in $A$. Let $a_{n}$ be
the infimum of all numbers $a\in(0,b_{n})$ for which the closed interval
$\left[a,b_{n}\right]$ is contained in $A$ and let $\widetilde{a}_{n}$
be the supremum of all numbers $a>b_{n}$ for which the closed interval
$\left[b_{n},a\right]$ is contained in $A$. Since $b_{n}^{\clubsuit}\notin A$
it follows that 
\begin{equation}
b_{n}^{\clubsuit}\le a_{n}.\label{eq:BclubAn}
\end{equation}
Since $A$ is open, it follows that $a_{n}\in\partial A$ and 
\begin{equation}
a_{n}<b_{n}<\widetilde{a}_{n}\le\infty\label{eq:shcw}
\end{equation}
and$\left(a_{n},\widetilde{a}_{n}\right)\subset A$. It may happen
that $\widetilde{a}_{n}=\infty$. Otherwise we have $\widetilde{a}_{n}\in\partial A$.
We see that all the hypotheses of Claim \ref{claim:Astuff} hold if
we choose $a=a_{n}$ and $\widetilde{a}=\widetilde{a}_{n}$. So that
claim enables us to assert that

\begin{equation}
\int_{a_{n}^{\clubsuit}}^{t}g(s)ds\le\int_{a_{n}^{\clubsuit}}^{t}f(s)ds\mbox{ \,\ for all\,}t\in\left[a_{n}^{\clubsuit},\max\left\{ \widetilde{a}_{n},a_{n}^{\clubsuit}+1\right\} \right],\label{eq:ubfe}
\end{equation}
and also 

\begin{equation}
g(a_{n}^{\clubsuit})\le f(a_{n}^{\clubsuit})\label{eq:gandflat}
\end{equation}
where we have defined $a_{n}^{\clubsuit}$ to be the largest integer
which satisfies $a_{n}^{\clubsuit}\le a_{n}$. Since $b_{n}^{\clubsuit}$
is an integer, it follows from (\ref{eq:BclubAn}) that 
\begin{equation}
b_{n}^{\clubsuit}\le a_{n}^{\clubsuit}.\label{eq:AnIsOk}
\end{equation}

We define $a_{n}^{\diamondsuit}$ to be $\infty$ if $\widetilde{a}_{n}=\infty$.
Otherwise we define $a_{n}^{\diamondsuit}$ to be the largest integer
which satisfies $a_{n}^{\diamondsuit}\le\max\left\{ \widetilde{a}_{n},a_{n}^{\clubsuit}+1\right\} $.
So of course 
\begin{equation}
a_{n}^{\clubsuit}+1\le a_{n}^{\diamondsuit}.\label{eq:flatdiamond}
\end{equation}
Then we define $A_{n}$ to be the non-empty interval $[a_{n}^{\clubsuit},a_{n}^{\diamondsuit})$.
Since $a_{n}^{\clubsuit}<a_{n}^{\diamondsuit}\le\max\left\{ \widetilde{a}_{n},a_{n}^{\clubsuit}+1\right\} $
we can immediately deduce from (\ref{eq:ubfe}) that the property
(\ref{eq:A-ineq}) holds for these particular choices of $n$ and
$A_{n}$.

Our definitions of $a_{n}^{\clubsuit}$ and $b_{n}^{\diamondsuit}$
are such that $a_{n}^{\clubsuit}\le a_{n}$ and $b_{n}\le b_{n}^{\diamondsuit}$.
So, in view of these inequalities and (\ref{eq:AnIsOk}) and (\ref{eq:shcw}),
we have that 
\[
b_{n}^{\clubsuit}\le a_{n}^{\clubsuit}\le b_{n}^{\diamondsuit}.
\]
This means that 
\begin{equation}
A_{n}\cup B_{n}=[b_{n}^{\clubsuit},\max\left\{ b_{n}^{\diamondsuit},a_{n}^{\diamondsuit}\right\} ).\label{eq:unionAnBn}
\end{equation}

We now have to consider the two possible ``subcases'' of the current
case, (i.e., the case where (\ref{eq:GivesO2}) does not hold). 

The first of these subcases is where $\widetilde{a}_{n}$ and therefore
also $a_{n}^{\diamondsuit}$ are both infinite. In this subcase, in
addition to our choices above of $B_{n}$ and $A_{n}$, we complete
Procedure P by choosing $\Omega_{n}$ and $\Gamma_{n}$ to both be
the empty set. We note that here we have obtained the outcome \textbf{(O-3)
}since (\ref{eq:unionAnBn}) now gives us that $A_{n}\cup B_{n}=[b{}_{n}^{\clubsuit},\infty)$.
In turn this obviously ensures that (\ref{eq:AllContained}) holds.
We have already checked that the sets $A_{n}$ and $B_{n}$ have property
(iv). It is also obvious that our sets have all the other properties
in pargraphs (i), (iii) and (v). Again, some of these are completely
trivial because $\Omega_{n}$ and $\Gamma_{n}$ are empty. 

Now we turn to the remaining subcase, namely where $\widetilde{a}_{n}$
and therefore also $a_{n}^{\diamondsuit}$are both finite. Recall
that above we applied Claim \ref{claim:Astuff}, choosing the quantities$a$
and $\widetilde{a}$ in its formulation to be $a_{n}$ and $\widetilde{a}_{n}$
respectively. Now, for these same choices of $a$ and $\widetilde{a}$,
since we now also have that $\widetilde{a}=\widetilde{a}_{n}$ is
finite and is in $\partial A$, we can apply the last part of Claim
\ref{claim:Astuff} to obtain useful information about the integer
$a_{n}^{\diamondsuit}$, whose definition above corresponds exactly
with the definition of $a^{\diamondsuit}$ in Claim \ref{claim:Astuff}.
Namely we can assert that at least one of the two integers $a_{n}^{\diamondsuit}$
and $a_{n}^{\diamondsuit}+1$ is an element of $B\setminus A$. So
we can define the new integer $b_{n+1}^{\clubsuit}$ by setting $b_{n+1}^{\clubsuit}=a_{n}^{\diamondsuit}+1$
if $a_{n}^{\diamondsuit}+1\in B\setminus A$ and otherwise setting
$b_{n+1}^{\clubsuit}=a_{n}^{\diamondsuit}$. 

In this subcase the sets $\Omega_{n}$ and $\Gamma_{n}$ will both
be non-empty. We choose $\Omega_{n}=[a_{n}^{\clubsuit},a_{n}^{\clubsuit}+1)$
and $\Gamma_{n}=[a_{n}^{\diamondsuit},a_{n}^{\diamondsuit}+1)$. This
completes Procedure P in this final subcase, and, as we will show
presently, its outcome is exactly as described by \textbf{(O-1).}

But, before showing that, we observe other consequences of Procedure
P in this subcase. Obviously, (\ref{eq:AllContained}) holds for our
four sets. Let us also now check that, as in previous (sub)cases,
these sets satisfy the properties appearing in paragraphs (i), (iii),
(iv) and (v) of Proposition \ref{prop:FourSequences} also this time.
This is obvious for the properties of (i). The property in (iii) is
exactly the same as (\ref{eq:flatdiamond}). We have already verified
above that $A_{n}$ and $B_{n}$ have the properties in (iv). Finally,
we deal with the three properties which are listed in (v). The first
two of them obviously hold, and the third property can be rewritten
as the inequality 
\begin{equation}
g(a_{n}^{\diamondsuit})\le f(a_{n}^{\clubsuit})\label{eq:Propv}
\end{equation}
which follows immediately from (\ref{eq:gandflat}) and (\ref{eq:flatdiamond}),
since $g$ is nonincreasing.

So now we turn to showing that, in this subcase, all the properties
described in the outcome \textbf{(O-1)} hold. 

Our definitions of course ensure that $A_{n}$, $B_{n}$, $\Omega_{n}$
and $\Gamma_{n}$ are each bounded and non-empty. The inequality (\ref{eq:bnClubGrows})
follows immediately from (\ref{eq:AnIsOk}) and (\ref{eq:flatdiamond})
combined with our definition of $b_{n+1}^{\clubsuit}$. 

So it only remains to prove (\ref{eq:ABGclub}). Our first step towards
doing this will be to show that the numbers $b_{n}^{\diamondsuit}$
and $\widetilde{a}_{n}$ generated by Procedure P in this subcase
satisfy

\begin{equation}
b_{n}^{\diamondsuit}\le\widetilde{a}_{n}.\label{eq:HearTilde}
\end{equation}
If $b_{n}$ is an integer then $b_{n}^{\diamondsuit}=b_{n}$ and (\ref{eq:HearTilde})
follows immediately from (\ref{eq:shcw}). If $b_{n}$ is not an integer
we will need the following longer explanation. Part of it is a sort
of analogue of our proof of (\ref{eq:gfaaN}) in Claim \ref{claim:Astuff}.
(Cf.~Remark \ref{rem:AboutStrongerVersion}.)

In this case $b_{n}^{\clubsuit}\le b_{n}<b_{n}^{\diamondsuit}$ and,
since $f^{q}-g^{q}$ takes the constant value $\left(f(b_{n}^{\diamondsuit}-1)\right)^{q}-\left(g(b_{n}^{\diamondsuit}-1)\right)^{q}$
on the interval $[b_{n}^{\diamondsuit}-1,b_{n}^{\diamondsuit})$ which
contains $b_{n}$ in its interior, we can substitute $t=b_{n}$ in
(\ref{eq:getheart}) and divide by $b_{n}^{\diamondsuit}-b_{n}$ in
order to show that $\left(f(b_{n}^{\diamondsuit}-1)\right)^{q}-\left(g(b_{n}^{\diamondsuit}-1)\right)^{q}\ge0$
which is of course the same as 
\begin{equation}
f(b_{n}^{\diamondsuit}-1)-g(b_{n}^{\diamondsuit}-1)\ge0.\label{eq:fgbheartz}
\end{equation}

Suppose that (\ref{eq:HearTilde}) does not hold, i.e., that $\widetilde{a}_{n}<b_{n}^{\diamondsuit}$.
Then, by (\ref{eq:shcw}), $\widetilde{a}_{n}$ must lie in the open
interval $\left(b_{n},b_{n}^{\diamondsuit}\right)$. We recall that
$\widetilde{a}_{n}$ is in $\partial A$ and is also the right endpoint
of the open interval $\left(a_{n},\widetilde{a}_{n}\right)$ which
is contained in $A$. This implies that the function $\varphi:[0,\infty)\to\mathbb{R}$,
which is defined by 
\begin{align*}
\varphi(t) & =\int_{0}^{t}f(s)ds-\int_{0}^{t}g(s)ds\,,
\end{align*}
must satisfy $\varphi(\widetilde{a}_{n})=0$ and $\varphi(t)>0$ for
all $t\in\left(a_{n},\widetilde{a}_{n}\right)$. So $\varphi$ cannot
be a nondecreasing function on any open interval which contains $\widetilde{a}_{n}$.
But this yields a contradiction, because at all points $t$ of the
open interval $\left(b_{n},b_{n}^{\diamondsuit}\right)$ and indeed
also for all points in the larger interval $[b_{n}^{\diamondsuit}-1,b_{n}^{\diamondsuit}]$
the function $\varphi$ satisfies 
\[
\varphi(t)=\int_{0}^{b_{n}^{\diamondsuit}-1}f(s)-g(s)ds+(t-b_{n}^{\diamondsuit})\left(f(b_{n}^{\diamondsuit}-1)-g(b_{n}^{\diamondsuit}-1)\right)
\]
and is in fact nonincreasing, in view of (\ref{eq:fgbheartz}). The
appearance of this contradiction proves that (\ref{eq:HearTilde})
also holds in this case, and therefore in all cases.

In view of (\ref{eq:HearTilde}), the integer $b_{n}^{\diamondsuit}$
also satisfies $b_{n}^{\diamondsuit}\le\max\left\{ \widetilde{a}_{n}.a_{n}^{\clubsuit}+1\right\} $
and so, since we have defined $a_{n}^{\diamondsuit}$to be the largest
integer which has this property, we must have 
\begin{equation}
b_{n}^{\diamondsuit}\le a_{n}^{\diamondsuit}.\label{eq:bhadim}
\end{equation}
This enables us to rewrite (\ref{eq:unionAnBn}) as the simpler formula
\[
A_{n}\cup B_{n}=[b_{n}^{\clubsuit},a_{n}^{\diamondsuit}).
\]
In view of our definitions of $\Omega_{n}$ (which is a subset of
$A_{n}$) and of $\Gamma_{n}$ and of $b_{n+1}^{\clubsuit}$, this,
in turn, shows that (\ref{eq:ABGclub}) holds. 

In view of (\ref{eq:bhadim}) and our definition of $b_{n+1}^{\clubsuit}$
we can also assert that 
\begin{equation}
b_{n}^{\diamondsuit}\le b_{n+1}^{\clubsuit}.\label{eq:bdiambclubPlusOne}
\end{equation}

Now that we have a complete explicit description of Procedure P for
all possible choices of the integer $b_{n}^{\clubsuit}$ in $B\setminus A$
and its ``label'' $n\in\mathbb{N}$, and now that we know that,
for each such $n$ and $b_{n}^{\clubsuit}$, the sets $A_{n}$, $B_{n}$,
$\Omega_{n}$ and $\Gamma_{n}$ have many, if not yet all, of the
properties listed in Proposition \ref{prop:FourSequences} we can
explicitly describe the iterative process which will complete the
proof of the proposition.

As the reader can guess, we first choose $n=1$ and let $b_{1}^{\clubsuit}$
be $0$ which is indeed a point in $B\setminus A$. We apply Procedure
P to obtain $A_{1}$, $B_{1}$, $\Omega_{1}$ and $\Gamma_{1}$. If
the outcome is \textbf{(O-1) }then we apply Procedure P to the point
$b_{2}^{\clubsuit}$ obtained from the previous step, of course now
using the ``label'' $n=2$. We proceed iteratively. After $n$ applications
of Procedure P, if the outcome is\textbf{ (O-1)} we apply the procedure
yet again to the point $b_{n+1}^{\clubsuit}$ which it previously
yielded, of course using the ``label'' $n+1$. If, after any number,
say $n_{0}$, of these iterations, the outcome is \textbf{(O-2)} or
\textbf{(O-3)} then we choose $A_{n}=B_{n}=\Omega_{n}=\Gamma_{n}=\emptyset$
for every $n>n_{0}$. Otherwise we continue indefinitely, thus obtaining
bounded non-empty sets $A_{n}$, $B_{n}$, $\Omega_{n}$ and $\Gamma_{n}$
for every $n\in\mathbb{N}$.

If for some particular $m\in\mathbb{N}$ the outcome of every one
of the first $m$ successive applications of Procedure P is\textbf{
(O-1)} then we obtain, from (\ref{eq:ABGclub}), that 
\begin{equation}
\bigcup_{n=1}^{m}\left(A_{n}\cup B_{n}\cup\Gamma_{n}\right)\supset\bigcup_{n=1}^{m}[b{}_{n}^{\clubsuit},b_{n+1}^{\clubsuit})=[0,b_{m+1}^{\clubsuit}).\label{eq:y6sv}
\end{equation}
If this is the case for every $m\in\mathbb{N}$ then (\ref{eq:bnClubGrows})
gives us that $\lim_{m\to\infty}b_{m}^{\clubsuit}=\infty$ and this
shows that the sequences of sets which we have generated satisfy (\ref{eq:UnionIsAll}).
Alternatively, in the situation already mentioned above, when, for
some $n_{0}$ the outcome of the $n_{0}$\textit{-th} iteration is
\textbf{(O-2)} or \textbf{(O-3) }then, with the help of (\ref{eq:y6sv})
and the properties of \textbf{(O-2)} or\textbf{ (O-3)} we see that 

\begin{align*}
\bigcup_{n=1}^{n_{0}}\left(A_{n}\cup B_{n}\cup\Gamma_{n}\right) & \supset A_{n_{0}}\cup B_{n_{0}}\cup\bigcup_{n=1}^{n_{0}-1}\left(A_{n}\cup B_{n}\cup\Gamma_{n}\right)\\
 & \supset[b_{n_{0}}^{\clubsuit},\infty)\cup[0,b_{n_{0}}^{\clubsuit})=[0,\infty).
\end{align*}
So (\ref{eq:UnionIsAll}) holds in all cases.

We still have to verify that the sets which we have constructed satisfy
the four properties listed in paragraph (ii) of the proposition. 

If the outcome of any one of the first $n$ successive applications
of Procedure P is not \textbf{(O-1)} then all these properties hold
trivially since then $A_{n+1}$, $B_{n+1}$, $\Omega_{n+1}$ and $\Gamma_{n+1}$
are all empty. So we only have to consider the situation where $B_{n}=[b_{n}^{\clubsuit},b_{n}^{\diamondsuit})$,
$A_{n}=[a_{n}^{\clubsuit},a_{n}^{\diamondsuit})$, $\Omega_{n}=[a_{n}^{\clubsuit},a_{n}^{\clubsuit}+1)$
and $\Gamma_{n}=[a_{n}^{\diamondsuit},a_{n}^{\diamondsuit}+1)$. First
we should keep in mind that all of the sets $A_{n+1}$, $B_{n+1}$,
$\Omega_{n+1}$ and $\Gamma_{n+1}$, whatever the outcome of Procedure
P for $n+1$, are either empty or contained in $[b_{n+1}^{\clubsuit},\infty)$.

If $b_{n+1}^{\clubsuit}=a_{n}^{\diamondsuit}+1$ then (by (\ref{eq:bhadim})
for $B_{n}$) all of $A_{n}$, $B_{n}$, $\Omega_{n}$ and $\Gamma_{n}$
are contained in $[b_{n}^{\clubsuit},b_{n+1}^{\clubsuit})$ and this
simultaneously establishes all four properties in this case. Otherwise
we have that $b_{n+1}^{\clubsuit}=a_{n}^{\diamondsuit}$ and then
we still have that $A_{n}$, $B_{n}$ and $\Omega_{n}$ are contained
in $[b_{n}^{\clubsuit},b_{n+1}^{\clubsuit})$. So it remains only
to show that $\Gamma_{n}<\Gamma_{n+1}$. If $\Gamma_{n+1}$ is not
empty then the outcome of Procedure P for $b_{n+1}^{\clubsuit}$ must
be \textbf{(O-1)} and $\Gamma_{n+1}$ must be the interval $[a_{n+1}^{\diamondsuit},a_{n+1}^{\diamondsuit}+1)$
and we can apply (\ref{eq:flatdiamond}) and (\ref{eq:AnIsOk}) both
with $n+1$ in place of $n$ to obtain that 
\begin{equation}
b_{n+1}^{\clubsuit}\le a_{n+1}^{\clubsuit}\le a_{n+1}^{\diamondsuit}-1.\label{eq:ste}
\end{equation}
In this case $\Gamma_{n}=[b_{n+1}^{\clubsuit},b_{n+1}^{\clubsuit}+1)$
and so (\ref{eq:ste}) suffices to establish that $\Gamma_{n}<\Gamma_{n+1}$.

We have now finished checking that our four sequences have all properties
stated in Proposition \ref{prop:FourSequences}, i.e., we have completed
the proof of that proposition.

We still have to use these sequences to establish (\ref{eq:toep}).

We shall do this with the help of six special non-negative functions,
namely $\varphi_{1}=f\chi_{\bigcup_{n\in\mathbb{N}}A_{n}}$, $\varphi_{2}=f\chi_{\bigcup_{n\in\mathbb{N}}B_{n}}$,
$\varphi_{3}=f\chi_{\bigcup_{n\in\mathbb{N}}\Omega_{n}}$, $\psi_{1}=g\chi_{\bigcup_{n\in\mathbb{N}}A_{n}}$,
$\psi_{2}=g\chi_{\bigcup_{n\in\mathbb{N}}B_{n}}$ and $\psi_{3}=g\chi_{\bigcup_{n\in\mathbb{N}}\Gamma_{n}}$. 

Let us recall (cf.~(\ref{eq:defgh}) and (\ref{eq:deff})) that the
function $f$ is related to our special sequence $x=\left\{ x_{n}\right\} _{n\in\mathbb{N}}$
in $E$ by the formula $f(n-1)=\left(1+\varepsilon\right)C(q)x_{n}^{*}$.
The following argument applies to each of the functions $\varphi_{j}$
for $j\in\left\{ 1,2,3\right\} $. The inequality $0\le\varphi_{j}\le f$
for all $t\in[0,\infty)$ implies that $\varphi_{j}^{*}\le f^{*}=f$
and so the nonincreasing rearrangement of the sequence $\left\{ \varphi_{j}(n-1)\right\} _{n\in\mathbb{N}}$,
which (cf.~Fact \ref{fact:FstarXstar}) is the same as the sequence
$\left\{ \varphi_{j}^{*}(n-1)\right\} _{n\in\mathbb{N}}$, satisfies
$\varphi_{j}^{*}(n-1)\le(1+\varepsilon)C(q)x_{n}^{*}$ for every $n\in\mathbb{N}$.
Therefore, by Claim \ref{claim:EisRI}, 
\begin{equation}
\left\{ \varphi_{j}(n-1)\right\} _{n\in\mathbb{N}}\in E\mbox{ and\,}\left\Vert \left\{ \varphi_{j}(n-1)\right\} _{n\in\mathbb{N}}\right\Vert _{E}\le(1+\varepsilon)C(q)C_{2}\left\Vert x\right\Vert _{E}\mbox{ for\,}j=1,2,3.\label{eq:PhiInE}
\end{equation}

In view of (\ref{eq:UnionIsAll}) the functions $\psi_{1}$, $\psi_{2}$
and $\psi_{3}$ satisfy
\begin{equation}
0\le g\le\psi_{1}+\psi_{2}+\psi_{3}.\label{eq:gAndPsis}
\end{equation}
 If we can show that, for $j=1,2,3$ the sequence $\left\{ \psi_{j}(n-1)\right\} _{n\in\mathbb{N}}$
is an element of $E$ then, by (\ref{eq:gAndPsis}) and Claim \ref{claim:Lattice},
we will be able to assert that the sequence $\left\{ g\left(n-1\right)\right\} _{n\in\mathbb{N}}$
is also an element of $E$ and 
\begin{equation}
\left\Vert \left\{ g(n-1)\right\} _{n\in\mathbb{N}}\right\Vert _{E}\le C_{2}\left\Vert \left\{ \sum_{j=1}^{3}\psi_{j}(n-1)\right\} _{n\in\mathbb{N}}\right\Vert _{E}\le C_{2}\sum_{j=1}^{3}\left\Vert \left\{ \psi_{j}(n-1)\right\} _{n\in\mathbb{N}}\right\Vert _{E}.\label{eq:rtohp}
\end{equation}
Recall from (\ref{eq:defgh}) that $\left\{ g(n-1)\right\} _{n\in\mathbb{N}}$
is the sequence $\left\{ y_{n}^{*}\right\} _{n\in\mathbb{N}}$ , which
is the nonincreasing rearrangement of our special sequence $y=\left\{ y_{n}\right\} _{n\in\mathbb{N}}$.
So Claim \ref{claim:EisRI} will enable us to deduce from (\ref{eq:rtohp})
that 
\begin{equation}
y\in E\mbox{ and\,}\left\Vert y\right\Vert _{E}\le C_{2}^{2}\sum_{j=1}^{3}\left\Vert \left\{ \psi_{j}(n-1)\right\} _{n\in\mathbb{N}}\right\Vert _{E}.\label{eq:rshcfh}
\end{equation}

So we have reduced the proof of our theorem to showing that $\left\{ \psi_{j}(n-1)\right\} _{n\in\mathbb{N}}\in E$
for $j\in\left\{ 1,2,3\right\} $ with appropriate norm estimates.
The first step for doing this is
\begin{claim}
\label{claim:PhisAndPsis}The nonincreasing rearrangements of the
functions $\varphi_{1}$, $\varphi_{2}$, $\varphi_{3}$, $\psi_{1}$,
$\psi_{2}$ and $\psi_{3}$ satisfy

\begin{equation}
\int_{0}^{t}\psi_{1}^{*}(s)ds\le\int_{0}^{t}\varphi_{1}^{*}(s)ds\mbox{ for all\,}t\ge0\label{eq:A-stuff-PhiPsi}
\end{equation}
 and also 

\begin{equation}
\int_{t}^{\infty}\left(\psi_{2}^{*}(s)\right)^{q}ds\le\int_{t}^{\infty}\left(\varphi_{2}^{*}(s)\right)^{q}ds\mbox{ for all\,}t\ge0\label{eq:B-stuff-PhiPsi}
\end{equation}

and also 
\begin{equation}
\psi_{3}^{*}(t)\le\varphi_{3}^{*}(t)\mbox{ for all\,}t\ge0.\label{eq:GammaOmegaStuff}
\end{equation}
\end{claim}
\begin{rem}
\label{rem:SimpleMentalPicture}The following proof of this claim
is mainly some quite straightforward applications of the properties
listed in paragraphs (ii) and (iv) of Proposition\ref{prop:FourSequences}.
In fact the following simple ``mental picture'' might make the said
proof rather more transparent or even perhaps quite unnecessary: Suppose
$w:[0,\infty)\to[0,\infty)$ is a nonincreasing right continuous function
and we create a new function $\omega:[0,\infty)\to[0,\infty)$ from
$w$ by replacing $w$ by $0$ on a finite or infinite sequence $\left\{ I_{n}\right\} $
of pairwise disjoint semi-open subintervals of $[0,\infty)$ and leaving
$\omega$ equal to $w$ everywhere else. Then we can obtain the graph
of $\omega^{*}$ from the graph of $\omega$ by starting from $0$
and pushing the first non-zero part of the graph of $\omega$ to the
left (if necessary) until it meets the $y$-axis, and then successively
``pushing'' every subsequent non-zero part of the graph of $\omega$
to the left until it meets the nearest non-zero part obtained by our
previous ``push''. In other words, we can think of the parts of
the graph lying above the intervals $I_{n}$ as being completely ``compressible''
and the rest of the graph as being ``incompressible'', and we squeeze
the whole graph to the left so that the parts above the $I_{n}$'s
simply disappear. 
\end{rem}
Despite the preceding remark, here is a detailed proof of Claim \ref{claim:PhisAndPsis}. 

Let $\mathbb{N}_{0}=\left\{ n\in\mathbb{N}:B_{n}\ne\emptyset\right\} $
and let $R=\sup$$\mathbb{N}_{0}$. For each $n\in\mathbb{N}_{0}$
let $\sigma_{n}:B_{n}\to[0,\infty)$ be a one to one measure preserving
map, in fact simply an affine map of the interval $B_{n}=[b_{n}^{\clubsuit},b_{n}^{\diamondsuit})$
onto an interval $B_{n}^{\prime}=[\gamma_{n},\gamma_{n}+\delta_{n})$
of the same length. More explicitly, we set $\delta_{n}=b_{n}^{\diamondsuit}-b_{n}^{\clubsuit}$
for each $n\in\mathbb{N}_{0}$, and then $\sigma_{1}(t)=t$ for all
$t\in B_{1}=[0,\delta_{1})=B_{1}^{\prime}$. And then, for each $n\ge2$
in $\mathbb{N}_{0}$, we set $\sigma_{n}(t)=t+\sum_{k=1}^{n-1}\delta_{k}-b_{n}^{\clubsuit}$,
so that $\sigma_{n}\left(B_{n}\right)=B_{n}^{\prime}=[\sum_{k=1}^{n-1}\delta_{k},\sum_{k=1}^{n}\delta_{k})$. 

Since the sets $B_{n}$ are pairwise disjoint (cf.~property (ii) of
Proposition \ref{prop:FourSequences} or (\ref{eq:bdiambclubPlusOne}))
we can combine the functions $\sigma_{n}$ to define a one to one
map $\sigma$ of the set $\bigcup_{n\in\mathbb{N}_{0}}B_{n}$ onto
the bounded or unbounded interval $[0,\sum_{n\in\mathbb{N}_{0}}\delta_{n})=\bigcup_{n\in\mathbb{N}_{0}}B_{n}^{\prime}$.
Whenever $t\in B_{n}$ for some $n\in\mathbb{N}_{0}$ we simply set
$\sigma(t)=\sigma_{n}(t)$. for each $t\in B_{n}$. Obviously $\sigma$
is also measure preserving. Furthermore, since $B_{n}<B_{n+1}$ (cf.~again
property (ii) of Proposition \ref{prop:FourSequences} or (\ref{eq:bdiambclubPlusOne})),
$\sigma$ is also a strictly increasing function on $\bigcup_{n\in\mathbb{N}_{0}}B_{n}$.
Its inverse $\sigma^{-1}$, which maps $[0,\sum_{n\in\mathbb{N}_{0}}\delta_{n})$
onto $\bigcup_{n\in\mathbb{N}_{0}}B_{n}$ is of course defined by
setting $\sigma^{-1}(t)=t-\sum_{k=1}^{n-1}\delta_{k}+b_{n}^{\clubsuit}$
whenever $t$ lies in $B_{n}^{\prime}$, and it is clearly right continuous
and measure preserving and strictly increasing. Since $f$ is nonincreasing
and right continuous, it follows that the composed function $f\circ\sigma^{-1}$
is nonincreasing and right continuous. Since $\sigma^{-1}$ only takes
values in $\bigcup_{n\in\mathbb{N}_{0}}B_{n}$, we also see that $f\circ\sigma^{-1}=\left(f\chi_{\bigcup_{n\in\mathbb{N}_{0}}B_{n}}\right)\circ\sigma^{-1}$.
Since $\sigma^{-1}$ is also a measure preserving map of $[0,\sum_{n\in\mathbb{N}_{0}}\delta_{n})$
onto $\bigcup_{n\in\mathbb{N}_{0}}B_{n}$, we see that $f\circ\sigma^{-1}$
has the same distribution function as the function $\varphi_{2}=$$f\chi_{\bigcup_{n\in\mathbb{N}_{0}}B_{n}}$.
All these properties ensure that 
\begin{equation}
f\circ\sigma^{-1}=\varphi_{2}^{*}\mbox{ on\,[0,\ensuremath{\sum_{n\in\mathbb{N}_{0}}\delta_{n}}).}\label{eq:fphisigma}
\end{equation}
If the sum $\sum_{n\in\mathbb{N}_{0}}\delta_{n}$ is finite, 
i.e., if the measure of $\bigcup_{n\in\mathbb{N}_{0}}B_{n}$ is finite,
then the set where $\varphi_{2}>0$ cannot be greater than $\sum_{n\in\mathbb{N}_{0}}\delta_{n}$
and so 
\begin{equation}
\mbox{If\,}\sum_{n\in\mathbb{N}_{0}}\delta_{n}<\infty\mbox{ then }\varphi_{2}^{*}(t)=0\mbox{\,\ for all\,}t\ge\sum_{n\in\mathbb{N}_{0}}\delta_{n}.\label{eq:Zero4BigT}
\end{equation}

For each $n\in\mathbb{N}_{0}$, in view of (\ref{eq:fphisigma}),
and an obvious affine change of variable 
\begin{eqnarray*}
\int_{B_{n}^{\prime}}\left(\varphi_{2}^{*}(s)\right)^{q}ds & = & \int_{B_{n}^{\prime}}\left(f(\sigma^{-1}(s))\right)^{q}ds\\
 & = & \int_{B_{n}^{\prime}}\left(f\left(s-\sum_{k=1}^{n-1}\delta_{k}+b_{n}^{\clubsuit}\right)\right)^{q}ds\\
 & = & \int_{b_{n}^{\clubsuit}}^{b_{n}^{\clubsuit}+\delta_{n}}\left(f(s)\right)^{q}ds=\int_{B_{n}}\left(f(s)\right)^{q}ds.
\end{eqnarray*}
Similarly, and a little more generally, if $t\in B_{n}^{\prime}$
for some $n\in\mathbb{N}_{0}$, then

\begin{eqnarray*}
\int_{B_{n}^{\prime}\cap[t.\infty)}\left(\varphi_{2}^{*}(s)\right)^{q}ds & = & \int_{[t,\sum_{k=1}^{n}\delta_{k})}\left(f(\sigma^{-1}(s))\right)^{q}ds\\
 & = & \int_{t}^{\sum_{k=1}^{n}\delta_{k}}\left(f\left(s-\sum_{k=1}^{n-1}\delta_{k}+b_{n}^{\clubsuit}\right)\right)^{q}ds\\
 & = & \int_{\sigma^{-1}(t)}^{b_{n}^{\clubsuit}+\delta_{n}}\left(f(s)\right)^{q}ds=\int_{B_{n}\cap[\sigma^{-1}(t),\infty)}\left(f(s)\right)^{q}ds.
\end{eqnarray*}
These previous calculations give us some useful formulae for $\Phi(t):=\int_{t}^{\infty}\left(\varphi_{2}^{*}(s)\right)^{q}ds$
for all $t\in[0,\infty)$. We first note that, by (\ref{eq:Zero4BigT}),
\begin{equation}
\Phi(t)=0\mbox{\,\ for all\,}t\ge\sum_{n\in\mathbb{N}_{0}}\delta_{n}\label{eq:BigPhi}
\end{equation}
if there are any such numbers $t$. Then, given any $t\in[0,\sum_{n\in\mathbb{N}_{0}}\delta_{n})$
there exists a unique $n(t)\in\mathbb{N}_{0}$ for which $t\in B_{n(t)}^{\prime}$.
So

\begin{eqnarray*}
\Phi(t) & = & \int_{t}^{\infty}\left(\varphi_{2}^{*}(s)\right)^{q}ds\\
 & = & \int_{t}^{\sum_{n\in\mathbb{N}_{0}}\delta_{n}}\left(\varphi_{2}^{*}(s)\right)^{q}ds\\
 & = & \int_{B_{n(t)}^{\prime}\cap[t,\infty)}\left(\varphi_{2}^{*}(s)\right)^{q}ds+\sum_{n=n(t)+1}^{R}\int_{B_{n}^{\prime}}\left(\varphi_{2}^{*}(s)\right)^{q}ds\\
 & = & \int_{B_{n(t)}\cap[\sigma^{-1}(t),\infty)}\left(f(s)\right)^{q}ds+\sum_{n=n(t)+1}^{R}\int_{B_{n}}\left(f(s)\right)^{q}ds.
\end{eqnarray*}
We can now repeat the preceding calculation exactly, but with $g$
in place of $f$ and therefore $\psi_{2}^{*}$ in place of $\varphi_{2}^{*}$.
It gives us that the function $\Psi(t):=\int_{t}^{\infty}\left(\psi_{2}^{*}(s)\right)^{q}ds$
satisfies

\begin{equation}
\Psi(t)=0\mbox{\,\ for all\,}t\ge\sum_{n\in\mathbb{N}_{0}}\delta_{n}\label{eq:BigPsi}
\end{equation}
if there are any such $t$, and, for all $t\in[0,\sum_{n\in\mathbb{N}_{0}}\delta_{n})$,
\[
\Psi(t)=\int_{B_{n(t)}\cap[\sigma^{-1}(t),\infty)}\left(g(s)\right)^{q}ds+\sum_{n=n(t)+1}^{R}\int_{B_{n}}\left(g(s)\right)^{q}ds.
\]
By (\ref{eq:B-ineq}) each integral in the sum which gives $\Psi(t)$
is dominated by the corresponding integral when $g$ is replaced by
$f$. This shows that $\Psi(t)\le\Phi(t)$ for all $t\in[0,\sum_{n\in\mathbb{N}_{0}}\delta_{n})$.
So, in view of (\ref{eq:BigPhi}) and (\ref{eq:BigPsi}) it holds
for $t\ge0$. This is exactly the inequality (\ref{eq:B-stuff-PhiPsi})
of Claim \ref{claim:PhisAndPsis}. 

The proof of (\ref{eq:A-stuff-PhiPsi}) is very similar to the preceding
proof of (\ref{eq:B-stuff-PhiPsi}). This time we use the index set
$\mathbb{N}_{1}=\left\{ n\in\mathbb{N}:A_{n}\ne\emptyset\right\} $
which may be slightly different from $\mathbb{N}_{0}$. We define
the sets $A_{n}^{\prime}$ for each $n\in\mathbb{N}_{1}$ by the same
formulae as we used for $B_{n}^{\prime}$ but this time with $\delta_{n}$
given by$\delta_{n}=a_{n}^{\diamondsuit}-a_{n}^{\clubsuit}$. (Note
that here we do not necessarily have $A_{1}^{\prime}=A_{1}$.) Again
we introduce a one to one piecewise affine strictly increasing measure
preserving map $\sigma$ of $\bigcup_{n\in\mathbb{N}_{1}}A_{n}$ onto
the interval $[0,\sum_{n\in\mathbb{N}_{1}}\delta_{n})=\bigcup_{n\in\mathbb{N}_{1}}A_{n}^{\prime}$
and use its inverse to find a formula for $\varphi_{1}^{*}$ which
is an analogue of (\ref{eq:fphisigma}). We then show that, for each
$n\in\mathbb{N}_{1}$, we have $\int_{A_{n}^{\prime}}\varphi_{1}^{*}(s)ds=\int_{A_{n}}f(s)ds$
and, more generally, $\int_{A_{n}^{\prime}\cap[0,t)}\varphi_{1}^{*}(s)ds=\int_{A_{n}\cap[0,\sigma^{-1}(t))}f(s)ds$
for all $t\in A_{n}^{\prime}$. We leave it to the reader to check
that the rest of the proof of (\ref{eq:A-stuff-PhiPsi}) proceeds
quite analogously to that of (\ref{eq:B-stuff-PhiPsi}), this time
with (\ref{eq:A-ineq}) playing a role analogous to that played by
(\ref{eq:B-ineq}) before.

For the proof of (\ref{eq:GammaOmegaStuff}) the relevant index set
is the set of all $n\in\mathbb{N}$ for which $A_{n}$ is non-empty
and bounded, and we denote it by$\mathbb{N}_{2}$ and its cardinality
by $R$. Here we introduce two different one-to-one piecewise affine
strictly increasing measure preserving maps $\sigma$ and $\widetilde{\sigma}$
which map $\bigcup_{n\in\mathbb{N}_{2}}\Omega_{n}$ or, respectively,
$\bigcup_{n\in\mathbb{N}_{2}}\Gamma_{n}$ onto $[0,R)$. More precisely,
$\sigma$ maps $[a_{n}^{\clubsuit},a^{\clubsuit}+1)$ onto $[n-1,n)$
and $\widetilde{\sigma}$ maps $[a_{n}^{\diamondsuit},a_{n}^{\diamondsuit}+1)$
onto $[n-1,n)$ for each $n\in\mathbb{N}_{2}$. Obviously $\varphi_{3}^{*}(t)=f(\sigma^{-1}(t))=\left(f\chi_{\bigcup_{n\in\mathbb{N}_{2}}\Omega_{2}}\right)(\sigma^{-1}(t))$
for all $t\in[0,R)$ and $\varphi_{3}^{*}(t)=0$ for any $t>R$. Similarly
$\psi_{3}^{*}(t)=g(\widetilde{\sigma}^{-1}(t))=\left(g\chi_{\bigcup_{n\in\mathbb{N}_{2}}\Gamma_{2}}\right)(\widetilde{\sigma}^{-1}(t))$
for all $t\in[0,R)$ and $\psi_{3}^{*}(t)=0$ for any $t>R$. The
inequality (\ref{eq:GammaOmegaStuff}) follows from these formulae
together with (\ref{eq:Propv}). 

This completes the proof of Claim \ref{claim:PhisAndPsis}.

\smallskip{}

It follows from (\ref{eq:A-stuff-PhiPsi}) and Theorem \ref{thm:Calderon}
that there exists $T\in\mathcal{L}_{1}\left((L^{1},L^{\infty})\right)$
for which $T\varphi_{1}=\psi_{1}$. Composing $T$ with the maps $P$
and $Q$ of Fact \ref{fact:PQ} we see that $QTP$ is a map in $\mathcal{L}_{1}\left((\ell^{1},\ell^{\infty})\right)$
which maps $\left\{ \varphi_{1}(n-1)\right\} _{n\in\mathbb{N}}$ to
$\left\{ \psi_{1}(n-1)\right\} _{n\in\mathbb{N}}$. In view of (\ref{eq:PhiInE})
and the fact that $E\in Int_{C_{2}}\left((\ell^{1},\ell^{\infty})\right)$
this shows that 
\begin{equation}
\left\{ \psi_{1}(n-1)\right\} _{n\in\mathbb{N}}\in E\mbox{ and\,}\left\Vert \left\{ \psi_{1}(n-1)\right\} _{n\in\mathbb{N}}\right\Vert _{E}\le(1+\varepsilon)C(q)C_{2}^{2}\left\Vert x\right\Vert _{E}.\label{eq:PsiOne}
\end{equation}

When we substitute $t=N-1$ for any $N\in\mathbb{N}$ in the inequality
(\ref{eq:B-stuff-PhiPsi}), it becomes 
\[
\sum_{n=N-1}^{\infty}\int_{n}^{n+1}\left(\psi_{2}^{*}(s)\right)^{q}ds\le\sum_{n=N-1}^{\infty}\int_{n}^{n+1}\left(\varphi_{2}^{*}(s)\right)^{q}ds
\]
 which is the same as $\sum_{n=N-1}^{\infty}\left(\psi_{2}^{*}(n)\right)^{q}\le\sum_{n=N-1}^{\infty}\left(\varphi_{2}^{*}(n)\right)^{q}$
or
\begin{equation}
\sum_{n=N}^{\infty}\left(\psi_{2}^{*}(n-1)\right)^{q}\le\sum_{n=N}^{\infty}\left(\varphi_{2}^{*}(n-1)\right)^{q}\mbox{ for every\,}N\in\mathbb{N}.\label{eq:StartPhiPsiTwo}
\end{equation}
Once again we keep in mind (Fact \ref{fact:FstarXstar}) that the
sequence $\left\{ \psi_{2}^{*}(n-1)\right\} _{n\in\mathbb{N}}$ is
the nonincreasing rearrangement of the sequence $\left\{ \psi_{2}(n-1)\right\} _{n\in\mathbb{N}}$
and similarly $\left\{ \varphi_{2}^{*}(n-1)\right\} _{n\in\mathbb{N}}$
is the nonincreasing rearrangement of $\left\{ \varphi{}_{2}(n-1)\right\} _{n\in\mathbb{N}}$.
So (\ref{eq:StartPhiPsiTwo}) together with (\ref{eq:PhiInE}) and
Claim \ref{claim:EasySqC} implies that 
\begin{equation}
\left\{ \psi_{2}(n-1)\right\} _{n\in\mathbb{N}}\in E\mbox{ and }\left\Vert \left\{ \psi_{2}(n-1)\right\} _{n\in\mathbb{N}}\right\Vert _{E}\le(1+\varepsilon)C(q)C_{1}C_{2}^{2}\left\Vert x\right\Vert _{E}.\label{eq:PsiTwo}
\end{equation}
It remains to use (\ref{eq:GammaOmegaStuff}) which, again keeping
Fact \ref{fact:FstarXstar} in mind, enables us to apply Claim \ref{claim:EisRI}
together with (\ref{eq:PhiInE}) to obtain that 
\begin{equation}
\left\{ \psi_{3}(n-1)\right\} _{n\in\mathbb{N}}\in E\mbox{ and }\left\Vert \left\{ \psi_{3}(n-1)\right\} _{n\in\mathbb{N}}\right\Vert _{E}\le(1+\varepsilon)C(q)C_{2}^{2}\left\Vert x\right\Vert _{E}.\label{eq:PsiThree}
\end{equation}
Our conclusions (\ref{eq:PsiOne}), (\ref{eq:PsiTwo}) and (\ref{eq:PsiThree})
are exactly the results that we said we would need in our discussion
immediately after (\ref{eq:gAndPsis}). Thus we have finally proved
that $y\in E$ and we can also substitute in (\ref{eq:rshcfh}) to
obtain that 
\begin{equation}
\left\Vert y\right\Vert _{E}\le C_{2}^{2}(1+\varepsilon)C(q)\left(C_{2}^{2}+C_{1}C_{2}^{2}+C_{2}^{2}\right)\left\Vert x\right\Vert _{E}.\label{eq:NearlyDone}
\end{equation}
Since our reasoning which establishes (\ref{eq:NearlyDone}) is valid
for every choice of $\varepsilon>0$, this shows that we have obtained
(\ref{eq:toep}) for the constant $C_{3}=C(q)C_{2}^{4}(C_{1}+2)$
and therefore completed the proof of Theorem \ref{thm:Main}. $\qed$

\section{\label{sec:Further}A consequence of the weak Fatou property and
Property $S_{q}(C)$}

As mentioned in the previous section, we turn here to elaborating
upon the remarks which were made there in the preamble to Theorem
\ref{thm:Main}.

We shall be using the set $\mathcal{M}$ of all linear maps $M:\ell^{\infty}\to\ell^{\infty}$
which are of the form 
\[
M\left(\left\{ x_{n}\right\} _{n\in\mathbb{N}}\right)=\left\{ \theta{}_{n}x_{\sigma(n)}\right\} _{n\in\mathbb{N}}
\]
where $\left\{ \theta_{n}\right\} _{n\in\mathbb{N}}$ is a sequence
which satisfies $\left|\theta_{n}\right|=1$ for all $n\in\mathbb{N}$
and $\sigma:\mathbb{N}\to\mathbb{N}$ is a one to one map of $\mathbb{N}$
onto itself. (We choose the notation $\mathcal{M}$ and $M$ here
in honour of Boris Mityagin since operators of forms similar to these
operators play an important role in his paper \cite{MityaginB1965}).
We observe that $\mathcal{M}$ obviously has the following properties:

(i) Each $M\in\mathcal{M}$ is an isometry on $\ell^{r}$ for every
$r\in[1,\infty]$, and thus, in particular, for $r=q$. 

(ii) For each sequence $x=\left\{ x_{n}\right\} _{n\in\mathbb{N}}\in\ell^{\infty}$
the nonincreasing rearrangements of $x$ and of $Mx$ are the same
sequence $x^{*}=\left\{ x_{n}^{*}\right\} _{n\in\mathbb{N}}$. 

(iii) Furthermore, if any such sequence $\left\{ x_{n}\right\} _{n\in\mathbb{N}}$
has finite support, or has no non-zero elements and satisfies $\lim_{n\to\infty}x_{n}=0$
then (cf.~Lemma~\ref{lem:EasyCalderon}) there exists an element
$M\in\mathcal{M}$ for which $x^{*}=Mx.$ 

In vew of the first two of these properties, the facts that $E$ has
property $S_{q}(C)$ and is contained in $\ell^{p}$ imply that
\begin{equation}
Mx\in E\mbox{ and\,}\left\Vert Mx\right\Vert _{E}\le C\left\Vert x\right\Vert _{E}\mbox{ for each\,}x\in E\mbox{ and\,}M\in\mathcal{M}.\label{eq:MEinE}
\end{equation}

\begin{lem}
\label{lem:LikeACLemma1}Let $\left\{ x_{n}\right\} _{n\in\mathbb{N}}$
be a sequence in $c_{0}$ and let $\left\{ y_{n}\right\} _{n\in\mathbb{N}}$
be a sequence with only finitely many non-zero elements. Suppose that
the non-increasing rearrangements of these two sequences satisfy 
\[
\sum_{n=1}^{m}y_{n}^{*}\le\sum_{n=1}^{m}x_{n}^{*}\mbox{ for every\,}m\in\mathbb{N}.
\]
Then there exists an linear map $T$ which is a convex combination
of elements of $\mathcal{M}$ such that $T\left(\left\{ x_{n}\right\} _{n\in\mathbb{N}}\right)=\left\{ y_{n}\right\} _{n\in\mathbb{N}}$.
\end{lem}
This result is well known. It can be deduced, for example, as a special
case of a slight modification of the argument used on pp. 226-227
of \cite{CalderonA1966} for the proof of Lemma 1 of \cite[p.~275]{CalderonA1966}.
It can also be easily deduced from two classical theorems of Birkhoff
and Hardy-Littlewood-P\'olya (see e.g., Theorem B and Theorem HLP
on pp.~233-234 of \cite{SparrG1978}. 
\begin{thm}
\label{thm:Enough}Let $q\in[1,\infty)$ and $C\ge1$ and $R\ge1$
be constants. Let $E$ be a normed sequence space which is contained
in $\ell^{q}$ and has property $S_{q}(C)$ and also property $WFP(R)$.
Then $E\in Int_{CR}\left((\ell^{1},\ell^{\infty})\right)$.
\end{thm}
\noindent \textit{Proof.} In view of Fact \ref{fac:KspaceIsInterp}
it will suffice if we show that $E$ is a $CR-K$ space with respect
to $\left(\ell^{1},\ell^{\infty}\right).$ (Of course this is also
known to be \textit{equivalent} to $E\in Int_{CR}\left((\ell^{1},\ell^{\infty})\right)$
but we do not need to use that fact here.) So let us suppose that
$x=\left\{ x_{n}\right\} _{n\in\mathbb{N}}$ is an arbitrary element
of $E$ and that the arbitrary sequence $y=\left\{ y_{n}\right\} _{n\in\mathbb{N}}$
satisfies 
\begin{equation}
K(t,y;\ell^{1},\ell^{\infty})\le K(t,x;\ell^{1},\ell^{\infty})\mbox{\,\ for all\,}t>0.\label{eq:KfunXY}
\end{equation}
We have to show that this implies that 
\begin{equation}
y\in E\mbox{ and\,}\left\Vert y\right\Vert _{E}\le CR\left\Vert x\right\Vert _{E}.\label{eq:shwnfh}
\end{equation}

If we apply (\ref{eq:KfunLoneLinf}) to the case where the measure
space is $\mathbb{N}$ equipped with counting measure, we obtain,
for each $\left\{ h_{n}\right\} _{n\in\mathbb{N}}\in\ell^{\infty}$,
that 
\begin{equation}
K(t,\left\{ h_{n}\right\} _{n\in\mathbb{N}};\ell^{1},\ell^{\infty})=\int_{0}^{t}h(s)ds\label{eq:Kfun4SeqSp}
\end{equation}
where $\left\{ h_{n}^{*}\right\} _{n\in\mathbb{N}}$ is the nonincreasing
arrangement of $\left\{ h_{n}\right\} _{n\in\mathbb{N}}$ and $h=\sum_{n=1}^{\infty}h_{n}^{*}\chi_{[n-1,n)}$.
(Here again we are using \ref{fact:FstarXstar}.) In particular, for
each $m\in\mathbb{N}$ this gives us that 
\begin{equation}
K(N,\left\{ h_{n}\right\} _{n\in\mathbb{N}};\ell^{1},\ell^{\infty})=\sum_{n=1}^{m}h_{n}^{*}.\label{eq:Kfun4N}
\end{equation}

For each $N\in\mathbb{N}$, the linear ``projection'' map $\Pi_{N}$
defined by $\Pi_{N}\left(\left\{ h_{n}\right\} _{n\in\mathbb{N}}\right)=\left\{ h_{n}^{(N)}\right\} _{n\in\mathbb{N}}$
is in $\mathcal{L}_{1}\left((\ell^{1},\ell^{\infty})\right)$. So,
by standard properties of $K$-functionals, we have 

\begin{equation}
K(t,\Pi_{N}y;\ell^{1},\ell^{\infty})\le K(t,y;\ell^{1},\ell^{\infty})\mbox{ for all\,}t>0.\label{eq:PiMakesSmaller}
\end{equation}

Let $\left\{ (y^{(N)})_{n}^{*}\right\} $ denote the nonincreasing
rearrangement of $\Pi_{N}y=\left\{ y_{n}^{(N)}\right\} _{n\in\mathbb{N}}$.
Then, in view of (\ref{eq:KfunXY}), (\ref{eq:PiMakesSmaller}) and
(\ref{eq:Kfun4N}), we have that 
\[
\sum_{n=1}^{m}(y^{(N)})_{n}^{*}\le\sum_{n=1}^{m}x_{n}^{*}\mbox{ for all\,}m\in\mathbb{N}.
\]
Since the sequence $\left\{ y_{n}^{(N)}\right\} _{n\in\mathbb{N}}$
has only finitely many non-zero elements, and since $x\in E\subset\ell^{q}\subset c_{0}$,
we can apply Lemma \ref{lem:LikeACLemma1} to obtain a linear map
$T$ which is a finite sum of the form $T=\sum_{j=1}^{J}\lambda_{j}M_{j}$.
where $M_{j}\in\mathcal{M},$ and $\lambda_{j}>0$ for all $j\in\left\{ 1,2,...,J\right\} $
and $\sum_{j=1}^{J}\lambda_{j}=1$, and which also satisfies $Tx=\left\{ y_{n}^{(N)}\right\} _{n\in\mathbb{N}}$.
Then we can use (\ref{eq:MEinE}) to deduce that $\left\{ y_{n}^{(N)}\right\} _{n\in\mathbb{N}}\in E$
and that 
\begin{eqnarray}
\left\Vert \left\{ y_{n}^{(N)}\right\} _{n\in\mathbb{N}}\right\Vert _{E} & = & \left\Vert \sum_{j=1}^{J}\lambda_{j}M_{j}x\right\Vert _{E}\nonumber \\
 & \le & \sum_{j=1}^{J}\lambda_{j}\left\Vert M_{j}x\right\Vert _{E}\le C\left\Vert x\right\Vert _{E}.\label{eq:NotForQuasinorms}
\end{eqnarray}
Since these inequalities hold for every $N\in\mathbb{N}$ and since
$E$ has property $WFP(R)$ we can deduce that $\left\{ y_{n}\right\} _{n\in\mathbb{N}}\in E$
and that $\left\Vert \left\{ y_{n}\right\} _{n\in\mathbb{N}}\right\Vert _{E}\le CR\left\Vert x\right\Vert _{E}$.
Thus we have established (\ref{eq:shwnfh}) and completed the proof
of the theorem. $\qed$
\begin{rem}
\label{rem:NotNecessary}As we briefly indicated in the preamble to
Theorem \ref{thm:Main}, that theorem is also applicable to some spaces
$E$ which do not have property $WFP(R)$. One example of such a space
is $\left(\ell^{p,\infty}\right)^{\circ}$ for $p\in(1,q)$, which
is defined to consist of all sequences $\left\{ x_{n}\right\} _{n\in\mathbb{N}}$
whose nonincreasing rearrangements satisfy 
\[
\lim_{m\to\infty}m^{1/p-1}\sum_{n=1}^{m}x_{n}^{*}=0
\]
 and is which is normed by 
\[
\left\Vert \left\{ x_{n}\right\} _{n\in\mathbb{N}}\right\Vert _{\ell^{p,\infty}}=\sup_{m\in\mathbb{N}}m^{1/p-1}\sum_{n=1}^{m}x_{n}^{*}.
\]
A simple calculation using the sequence $\left\{ n^{-1/p}\right\} _{n\in\mathbb{N}}$
and its truncations shows that $\left(\ell^{p,\infty}\right)^{\circ}$
does not have property $WFP(R)$ for any $R>0$. But it is obviously
a $1-K$ space with respect to the couple $\left(\ell^{1},\ell^{\infty}\right)$
and therefore in $Int_{1}\left((\ell^{1},\ell^{\infty})\right)$. 

It must be conceded that Theorem \ref{thm:Main} is not of any real
use for dealing with this particular space because apparently the
easiest way to see that $\left(\ell^{p,\infty}\right)^{\circ}$ has
property $S_{q}(C)$ for some $C$ is to observe that this space is
the closure of $\ell^{1}$ in the Weak $\ell^{p}$ space $\ell^{p,\infty}$
which, by the reiteration theorem for Lions-Peetre spaces is, to within
equivalence of norms, equal to $\left(\ell^{1},\ell^{q}\right)_{\theta,\infty}$
where $\theta$ satisfies $1/p=(1-\theta)+\theta/q$. This shows that
$\left(\ell^{p,\infty}\right)^{\circ}\in Int_{C}\left((\ell^{1},\ell^{q})\right)$
for some $C>0$ and so has property $S_{q}(C)$ by Theorem \ref{thm:Easy}. 
\end{rem}

\section{\label{sec:Appendices}Appendices}

\subsection{\label{sec:MarcelOlof}Remarks about the classical interpolation
theorems of Marcel Riesz and Riesz-Thorin}

The readers of this paper are surely very familiar with the terminology
and statement and proof of the Riesz-Thorin interpolation theorem
\cite{ThorinG1939} concerning linear operators $T$ which act on
spaces of complex valued functions and are bounded from $L^{p_{0}}$ to
$L^{q_{0}}$ and also from $L^{p_{1}}$ to $L^{q_{1}}$ for given exponents
$p_{0}$, $p_{1}$, $q_{0}$, $q_{1}$ in $\left[1,\infty\right]$
and are shown to consequently be bounded from $L^{p_{\theta}}$ into
$L^{q_{\theta}}$ and satisfy the norm estimates 
\begin{equation}
\left\Vert T\right\Vert _{L^{p_{\theta}}\to L^{q_{\theta}}}\le\left\Vert T\right\Vert _{L^{p_{0}}\to L^{q_{0}}}^{1-\theta}\left\Vert T\right\Vert _{L^{p_{1}}\to L^{q_{1}}}^{\theta}\label{eq:LogConvex}
\end{equation}

But some readers may perhaps be a little less familiar with Marcel
Riesz's earlier version \cite{RieszM1926} of this theorem, which
makes an exactly analogous statement about operators acting on $L^{p}$ spaces
of \textit{real} valued functions, but subject to the additional requirements
\begin{equation}
p_{0}\le q_{0}\mbox{\,\ and\,}p_{1}\le q_{1}\label{eq:AdReq}
\end{equation}
on the above mentioned exponents. Various equivalent statements of
it, and their proofs, can be seen in \cite[Section 1.4, pp.~13--16]{BrudnyiYKrugljakN1991}
and (more concisely) in \cite[pp.~73-74]{PeetreJ1979dm} as well as,
of course, in \cite[pp.~466--472]{RieszM1926}). It can very easily
be deduced from the Riesz-Thorin theorem, but only with a weaker variant
of the estimate (\ref{eq:LogConvex}). 

Also there is a counterexample due to Marcel Riesz (see \cite[Section 39, pp.~495--496]{RieszM1926}
or \cite[Remark 1.4.5, pp.~16--17]{BrudnyiYKrugljakN1991}) which
shows that, for certain choices of exponents which do not satisfy
(\ref{eq:AdReq}), the conclusion of his theorem, in particular the
estimate (\ref{eq:LogConvex}), does not always hold.

We wish to formulate our results in this paper for various spaces
of sequences or functions over both the complex and real fields. Fortunately
whenever we need to consider interpolation of operators acting on
real or complex $L^{p}$ spaces or $\ell^{p}$ spaces the relevant
exponents will always satisfy (\ref{eq:AdReq}) (even with equality),
and so, regardless of whether our spaces are over the real field or
the complex field, we will be able to ensure, via one or the other
of the two above mentioned interpolation theorems, that we obtain
the expected conclusion, including exactly the norm estimate (\ref{eq:AdReq}).

\subsection{\label{subsec:DiffDefnInterpSp}
Remarks related to our slightly different
definitions of interpolation spaces and $K$-spaces}

The reader may choose (and maybe we should have made the same choice)
to simply confine attention to those normed spaces $A$ considered
in Definition \ref{def:C-interp} and Definition \ref{def:Kspace}
which also happen to be Banach spaces and for which the inclusions
$A_{0}\cap A_{1}\subset A\subset A_{0}+A_{1}$ are continuous. As
already pointed out in Remark \ref{rem:VeryImportant}, our results
(in which the relevant space may be denoted by $E$ rather than $A$)
apply, without any need of any modification of our proofs, to the
special cases of these spaces, which are the kinds of interpolation
spaces and $K$ spaces which are usually considered in the literature.
We are certainly not the first nor the only authors to consider these
kinds of less stringent definitions of interpolation spaces. See for
example the work of Paul Kr\'ee \cite{KreeP1968}. 

It is easy to show that any space $A$ which is a normed interpolation
space with respect to some Banach couple $\left(A_{0},A_{1}\right)$
(as specified in Definition \ref{def:C-interp}) must contain $A_{0}\cap A_{1}$
and must be contained in $A_{0}+A_{1}$. But exotic examples show
that if $A$ is merely a normed interpolation space or a normed $K$-space,
these inclusions are not necessarily continuous. (Take $A$ to be
$A_{0}\cap A_{1}$ or $A_{0}+A_{1}$ or $\left(A_{0},A_{1}\right)_{\theta,p}$
and use a Hamel basis of $A$ to define a norm on $A$ which is not
equivalent to its usual norm.) 

If $A$ is a normed $C-K$ space with respect to a Banach couple $\left(A_{0},A_{1}\right)$
and some positive constant $C$ then it is not difficult to deduce
that the inclusions $A_{0}\cap A_{1}\subset A\subset A_{0}+A_{1}$
must both be continuous. The argument uses the inequalities 
\[
K(t,x;A_{0},A_{1})\le\min\left\{ 1,t\right\} \le K(t,y;A_{0},A_{1})
\]
which hold for every $t>0$ whenever $x$ is a norm $1$ element 
of $A_{0}\cap A_{1}$ and
$y$ is a norm 1 element of $A_{0}+A_{1}$. 

There are obvious examples of spaces $A$ which satisfy the conditions
stated in Definition \ref{def:C-interp} or in Definition \ref{def:Kspace}
but are not complete. For example the space $\ell^{1}$ equipped with
the norm of $\ell^{p}$ for some $p\in(1,\infty]$ is a normed $1$-interpolation
space (obviously, by the classical results recalled in Subsection\ref{sec:MarcelOlof})
and therefore, by Theorem\ref{thm:CalderonExplicit}, is also a $1-K$
space. But it is not complete. We can obtain similar examples, with
respect to suitable not too ``trivial'' Banach couples $\left(A_{0},A_{1}\right)$,
by choosing a Banach space $B$ which satisfies one or more of the
conditions in these definitions and choosing $A$ to be $B\cap A_{0}\cap A_{1}$
equipped with the norm of $B$. Or we can choose $A=B\cap E$, again
with the norm of $B$, where $E$ is some other suitably chosen interpolation
or $K$ space with respect to $\left(A_{0},A_{1}\right)$.

It would be interesting to discuss further connections and further
differences between the spaces defined by Definitions \ref{def:C-interp}
and \ref{def:Kspace} and those defined by the usual versions of these
definitions, and to consider other relevant examples. It should be
noted that in some cases an interpolation space or a $K$ space is
necessarily also, respectively, a $C$-interpolation space or a $C-K$
space for some constant $C$. But these matters are not a priority
of this paper. Material relevant to such a discussion can be found,
e.g., in various parts of \cite{AronszajnNGagliardoE1965}, including
in particular, Theorem 6.XI on p.~73. See also Theorem 6.1 on pp.~70-71
of \cite{CwikelMNilssonP2003} concerning complete $K$ spaces with
respect to couples of Banach lattices.

\end{document}